\documentstyle[amstex,amssymb,amscd,12pt]{amsart}
\newtheorem {theorem} {Theorem} [section]
\newtheorem {cor} [theorem] {Corollary}
\newtheorem {defn} {Definition}
\newtheorem {prop}
[theorem] {Proposition}
\newtheorem {lemma} [theorem] {Lemma}

\newtheorem {conj} [theorem] {Conjecture}

\newcommand{\codim}{\operatorname{codim}}

\newcommand{\spax}{\operatorname{span}}
\newcommand{\proj}{{\Bbb{P}}}
\newcommand{\aff}{{\Bbb{A}}}

\newcommand{\nats}{{\Bbb{N}}}
\newcommand{\ints}{{\Bbb{Z}}}
\newcommand{\f}{{\Bbb{F}}}
\newcommand{\cchar}{\operatorname{char}}
\newcommand{\z}{{\cal Z }}
\title [Geometric Conjectures on the Hilbert function of fat points ]
 {The geometric interpretation of Fr\"oberg-Iarrobino conjectures
on infinitesimal neighbourhoods of points in projective space}
\author {KAREN A. CHANDLER}
\begin{document}
\maketitle

\begin{abstract}
The study of infinitesimal deformations of a variety embedded in projective space requires, at ground level,
that of deformation of a collection of points, as specified by a zero-dimensional scheme.
Further, basic problems in infinitesimal interpolation correspond directly to the analysis of such schemes.

An optimal Hilbert function of a collection of infinitesimal neighbourhoods of points in projective space is
suggested by algebraic conjectures of R. Fr\"oberg and A. Iarrobino.
We discuss these conjectures from a geometric point of view.

The conjectures give, for each such collection, a function
(based on dimension, number of points, and order of each neighbourhood)
which should serve as an upper bound to its Hilbert function (Weak Conjecture).
The Strong Conjecture 
 predicts when the upper bound is sharp, in the case  of 
equal order throughout.  
In general we refer to the equality of the Hilbert function 
of a collection of infinitesimal neighbourhoods with that of the corresponding conjectural function
as the Strong Hypothesis.

We interpret  these conjectures and hypotheses as accounting for the infinitesimal neighbourhoods
of projective subspaces naturally occurring in the base locus of a linear system with prescribed
singularities at fixed points. 
We develop techniques and insight toward
the conjectures' verification and refinement.

The main result 
gives an  an upper bound on the Hilbert 
function of a collection of infinitesimal neighbourhoods in~$\proj^n$ based on Hilbert functions of 
certain such subschemes of~$\proj^{n-1}$.
Further, equality occurs exactly when the scheme has only the expected linear obstructions to the
linear system at hand.
It follows 
 that an infinitesimal neighbourhood scheme  obeys the Weak Conjecture provided that
the schemes identified in codimension~one satisfy  the Strong Conjecture.

This observation is then applied 
to show that 
the Weak Conjecture does hold valid in~$\proj^n$ for~$n \leq 3$.
The main feature here is that the result is obtained although the Strong Hypothesis 
is not known to hold generally in~$\proj^2$ and, further, $\proj^2$ presents
special exceptional cases.
Consequences of the main result  in higher dimension are then examined. 
 We note, then, that the full weight of the Strong Conjecture (and validity of the
Strong Hypothesis) are not necessary toward using the main theorem in the next dimension.

Further, we exhibit general situations in which the Strong Hypothesis
does hold valid: when the sum of orders of vanishing is not too large 
compared with degree; and in the case of equal multiplicities~$k$ in
degree~$k+1$.  

On the other hand, we construct (classes of) counterexamples to the Strong Conjecture.

We end with the observation 
 of how our viewpoint on the Strong 
Hypothesis
pertains to extra algebraic information: namely, on the structure
of the minimal free resolution of an ideal generated by linear forms.
\end{abstract}

\section{Introduction.}
Let~${\cal K}$ be an infinite field, and~$\proj^n = \proj^n_{\cal K}$.

The Hilbert function of a scheme~${\cal Z}$ embedded in projective
space evaluates in each degree the codimension of the graded 
piece of the ideal of~${\cal Z}$ with respect to the relevant homogeneous co\"ordinate ring.
One seeks in general to establish how the information provided by the
Hilbert function describes the geometry of the scheme and its embedding.

An {\bf infinitesimal neighbourhood} of a variety~$X$ with respect to an embedding in~$\proj^n$
is a scheme defined by a power~${\cal I}_X^k$ of its ideal sheaf.
We shall also refer to such a scheme as a (full) multiple subvariety of~$\proj^n$
of (overall) multiplicity~$k$, and as a ``fat subvariety'' of~$\proj^n$.
Each infinitesimal neighbourhood of a variety then refers to the extent of
singularity of a hypersurface through the variety itself. 
We study here the Hilbert function of a collection of infinitesimal neighbourhoods of
points in projective space.

An immediate motivation for this investigation is given by
{\bf infinitesimal interpolation}.
The Hilbert  function of a collection of infinitesimal neighbourhoods of points measures the number of
linear conditions imposed on the linear system of hypersurfaces of each
degree given by the requirement to vanish to specified order at each point
of a generic subset.
These data tell (in appropriate characteristic, say) the extent to which
it is possible to interpolate the values of a polynomial of given degree, together
with its partial derivatives up to specified orders, to a collection
of points in affine space. (See \cite{cil}, \cite{gs}, for example.)
(For ``inappropropriate characteristic'' the same principle applies, subject
to a modification of the notion of derivative; see \cite{ik}.)

Moreover, the study of such Hilbert functions is a basic starting point
in that of an {\bf infinitesimal deformation} of a (higher-dimensional!)
variety~$X$ embedded in a projective space.  (See, e.g.  \cite{me0}.)
For example, to estimate (or evaluate)  cohomologies of twists of
the $(k-1)$th symmetric power of the conormal bundle of~$X$
in the projective space one may examine those of the scheme
defined by~${\cal I}_X^k$.
The standard method of hyperplane slicing gives cohomological
data on this scheme from those of a lower-dimensional one.
However, in low degree of twisting (the most interesting!) the
standard approach is far too crude.  
We shall focus here on this phenomenon
 and
work toward refining the technique (see also \cite{me0}, \cite{me}, \cite{me2}, \cite{me3}, \cite{me4}, and \cite{mernc}),
and thereby advance the theory of infinitesimal deformation.

Special attention is paid here to  the situation of a generic collection of multiple
points (i.e., the support is a generic subset of projective space).
Some motivation for this restriction is evident: we do have the 
conjectures defined below as guidance.  Surely the problem of finding the
Hilbert function is made easier by having the freedom to choose 
generic points, and any upper bound obtained on the function 
for generic points gives automatically a bound for an arbitrary
collection of points.  
But we also proceed here with an eye 
toward developing tools applicable to the Hilbert function of any collection of
multiple points (or multiple varieties); such as identifying which such schemes 
have the maximal possible Hilbert function.
  For example we find
in~\cite{mernc} that a variation of the technique introduced here
applies well to a collection of multiple points lying on a rational
normal curve (which ought, according  to conjectures of Catalisano and Gimigliano \cite{cg}, to give the ``worst''
Hilbert function amongst sets of points in linearly general position).

A zero-dimensional subscheme~${\cal Z} \subset \proj^n$ is said 
to have the {\bf maximal rank property} if its Hilbert function
is as simple as possible: for each degree~$m$, 
either~${\cal Z}$ does not lie on an~$m$-ic hypersurface or~${\cal Z}$ 
imposes~$\deg {\cal Z}$ (i.e. {\bf independent}) conditions on the linear 
system of~$m$-ics.
One is led to consider, then, which~${\cal Z}$ 
{\it do not} enjoy this property?
One expects, at least conjecturally, 
that for a generic such scheme~${\cal Z}$ the failure of maximal rank in a given 
degree should occur when the base locus of the linear system of 
of~$m$-ics through~${\cal Z}$ is forced to contain a positive dimensional 
scheme whose intersection with~${\cal Z}$ itself cannot impose the ``expected 
number'' of conditions on the system. 

For a generic collection~${\cal Z}$ of multiple points in
projective space,  the inductive procedure (m\'ethode d'Horace
diff\'erentielle) of J. Alexander and A. Hirschowitz allows
one to deduce maximal rank in a given degree from  maximal
rank conditions in lower degree and in lower dimension.
This idea is used in \cite{ah5} to 
obtain asymptotic results on maximal rank.

But in low degree~$m$ (compared to order of vanishing) 
such a scheme~${\cal Z}$ cannot
impose independent conditions on~$m$-ics, due to visible linear obstructions.  
Specifically, suppose that~${\cal Z}$ contains two points, 
of multiplicities~$j, k$, and 
take~$m \leq j+k-2$.   The line~$L$ between the two points meets~${\cal Z}$
in a subscheme of degree~$j+k > m+1$ which then cannot
impose independent conditions on~$m$-ics (i.e., $L$ itself
imposes only~$m+1$ conditions) and hence neither does~${\cal Z}$.

Therefore, a key issue on obtaining information on the Hilbert function of
such a scheme, such as finding explicit (better yet, sharp) conditions for
maximal rank, is to study cases in which maximal rank is obstructed by linear subspaces
spanned by subsets of the set of points in the scheme.

Conjectures of R. Fr\"oberg and A. Iarrobino give a proposed
value (Strong Conjecture) or upper bound (Weak Conjecture) for such
a Hilbert function.
These conjectures arise indirectly from an algebraic conjecture of Fr\"oberg \cite{f}.
He studies an ideal generated by a generic collection of forms, and 
asserts that its behaviour may be quantified (or at least estimated)
by a natural generalisation of the formula for complete intersection ideals.
Iarrobino further asserts \cite{iar} that,
 up to an identifiable  region of cases,
the conjecture of Fr\"oberg should apply to  an ideal 
generated by a generic collection 
of powers of linear forms of equal degree.
(Of course, one must certainly exclude situations such as $p$th
powers in characteristic~$p$, by virtue of the ``Freshman's Dream Theorem''!)
An application of Macaulay duality, given by Emsalem and Iarrobino, to
the case of generic powers of linear forms yields the conjectures
on multiple points \cite{ei}.

We present here a direct geometric interpretation of these conjectures.
Namely, the conjectural Hilbert function of multiple points
reflects circumstances under which
 (multiple) planes spanned by subsets must appear in
the base locus at issue, according to Lagrange-Hermite, say.
In particular, we argue that the Strong Conjecture for multiple 
points  corresponds to situations in which the {\it only} obstructions
to the scheme's imposing independent conditions are the ``obvious
linear ones'', and that the Weak Conjecture amounts to {\it counting} such
linear obstructions.

\medskip

Let us recall:

\begin{defn}
For a subscheme~$\z \subset \proj^n$ with ideal sheaf~${\cal I}_{\z}$
the {\bf Hilbert function} of~$\z$ (as a function of~$m$) is given
by
$$h_{\proj^n}(\z,m):=
\dim H^0(\proj^n, {\mathcal O}_{\proj^n}(m))-\dim H^0(\proj^n,{\mathcal I}_\z(m)). $$
\end{defn}

\begin{defn}
Take a variety~$X \subset \proj^n$ and~$k \in \nats$.  The {\bf $(k-1)$th infinitesimal neighbourhood } of~$X$
(with respect to~$\proj^n$)
 is the scheme given by~${\cal I}^k$, where~${\cal I}$
is the ideal sheaf of~$X$.
We shall denote this scheme by~$X^k \subset \proj^n$.
\end{defn}

So, for example, $X^0 = \emptyset$, and~$X^1=X$.

Note that for~$p \in \proj^n$ the degree of~$\{p\}^k$ is~${{n+k-1} \choose n}$.
Hence, for an~$r$-dimensional variety~$X \subset \proj^n$, $\deg X^k = {{n+k-r-1} \choose {n-r}} \deg X$.

For brevity, when the ambient projective space of embedding is clear, we shall refer to~${X^k \subset \proj^n}$
as~$X^k$, a~$k$-uple subscheme (or a ``fat variety'' of multiplicity~$k$). 

\begin{defn}\label{ascheme}
Given~$A=(k_1, \ldots, k_d) \in \nats^d$ 
define an~$A$-{\bf subscheme of}~$\proj^n$
as a union of~$\{p_1\}^{k_1} \cup  \ldots \cup \{ p_d\}^{k_d}$
where~$\{p_1, \ldots, p_d\}$ is a set of~$d$ points in~$\proj^n$.
We shall say that an~$A$-scheme is {\bf homogeneous} (respectively, {\bf quasihomogeneous})
if~$k_1=\ldots=k_d$ (respectively, after perhaps reordering, $k_2= \ldots=k_d$).
(For emphasis, we may refer to a scheme as having {\bf mixed multiplicities} if it is not necessarily homogeneous.)
\end{defn}

The Fr\"oberg-Iarrobino Conjectures (see Section~4) refer to a function~$G(d,A,n+1)_m$
and its correspondence with the Hilbert function of an $A$-subscheme of~$\proj^n$
(supported on~$d$ points).
The Weak Conjecture (Conjecture~\ref{wpt}) asserts that for each~$d$-uple~$A \in \nats^d$ 
and each~$A$-subscheme~${\cal Z} \subset \proj^n$ we have
$$h_{\proj^n}( {\cal Z}, m ) \leq G(d,A,n+1)_m$$
for each degree~$m$; while the Strong Conjecture (Conjecture~\ref{spt})
gives numerical conditions under which equality is predicted to hold
between the two functions in the case of homogeneous schemes with generic support.
Here we provide techniques for analysing the Hilbert function of a collection of multiple
points and compare with the properties of the proposed function~$G$.

The Strong Fr\"oberg-Iarrobino Conjecture deals only with homogeneous schemes.
We aim, further, to find the Hilbert function of a scheme of mixed multiplicities.
Moreover (and unfortunately)  we shall exhibit (infinitely many) counterexamples to the Strong Conjecture below (Section~10);
that is, generic homogeneous  subschemes of~$\proj^n$ (for~$n=4, 5, 6$) whose Hilbert functions do not agree
with the prescribed values.
One seeks to: refine the Strong Fr\"oberg-Iarrobino Conjecture on homogeneous
schemes, and then to extend to those of mixed multiplicities.
According to the evidence presented here, this requires the identification of nonlinear positive dimensional varieties lying
in the base locus of a linear system of hypersurfaces through a multiple point scheme.
In general, we refer to the {\bf Strong Fr\"oberg-Iarrobino Hypothesis} (on a given case) 
as the supposition 
that a generic~$A$-subscheme of~$\proj^n$ has Hilbert function that agrees with 
the corresponding
Fr\"oberg-Iarrobino function (in a given degree). 
We shall also say  that the Strong Hypothesis {\it applies} to a given~$A$-scheme provided 
that its Hilbert function is equal to the conjectured value.

We obtain the following (Section~6): 

\begin{theorem} \label{siwa} Let~$n, m, d \in \nats$ and~$A=(k_1, \ldots, k_d) \in \nats^d$.
Take~$C_{ji}=((k_1+i-m)^+, \ldots,( k_j+i-m)^+)$ for each~$j=1, \ldots, d-1$, and~$i=0, \ldots, k_j-1$.

Suppose that the Strong Fr\"oberg-Iarrobino Hypothesis in~$\proj^{n-1}$ is verified by each generic~$C_{ji}$-subscheme of~$\proj^{n-1}$
in degree~$i$, for~$j=1, \ldots, d-1$ and~$i=0, \ldots, k_j-1$.

Then:\\
a) Each~$A$-subscheme of~$\proj^n$ satisfies the Weak Conjecture in degree~$m$, and \\
b) A generic~$A$-subscheme of~$\proj^n$ satisfies the Strong Hypothesis in degree~$m$ if and only if
it displays only the expected linear obstructions in degree~$m$ (See Definition~\ref{exp}) \\
c) If~${\cal Z} \subset \proj^n$ is any $A$-subscheme of~$\proj^n$
then~${\cal Z}$ verifies the Strong Hypothesis in degree~$m$ provided that
it admits only the expected linear obstructions in degree~$m$.
\end{theorem}

We shall observe in Proposition~\ref{siww} that we may ``homogenise'' this result to deal with the Strong
Conjecture itself.

The notion of expected linear obstructions, given in Definition~\ref{exp} (Section~5), is 
 based simply on the prediction of B\'ezout on how a line 
(and whence multiple lines as well as higher dimensional linear subspaces) must
appear in the base locus of a linear system if its intersection with
that base locus has sufficient degree.

Theorem~\ref{siwa} implies that
 one may determine whether a given $A$-subscheme of~$\proj^n$  satisfies the Weak Conjecture
by finding analogous subschemes of~$\proj^{n-1}$ for which the Strong Hypothesis applies. 
In particular, to verify the Weak Conjecture in~$\proj^{n}$, it suffices that the Strong Hypothesis
applies to sufficiently many (and identifiable) subschemes of~$\proj^{n-1}$.

For example, in~$\proj^2$ the Strong Conjecture displays ``extra exceptions'' to the expected maximal
rank of a (homogeneous) multiple subscheme. 
 These are extended in~\cite{cm1} to conjectures on subschemes of~$\proj^2$ of mixed multiplicities.
Progress on these conjectures on~$\proj^2$ has been made recently (see~Section~3), but the main problem remains open,
even for homogeneous schemes.
Nevertheless, we employ  Theorem~\ref{siwa} 
to obtain (Section~7):

\begin{theorem}\label{p3}
The Weak Fr\"oberg-Iarrobino conjecture holds valid
 in $\proj^n$ for~$n \leq 3$.
\end{theorem}

Hence we see that the ``full strength'' of the Strong Hypothesis is not necessary toward 
verifying the Weak Conjecture in the next dimension from Theorem~\ref{siwa}.


\medskip

The main tool in finding  an upper bound for the Hilbert 
function of a given collection of infinitesimal neighbourhoods (homogeneous or otherwise)
is presented in Lemma~\ref{key}.
 Here an inductive strategy on comparing the Hilbert function of a scheme, say~${\cal Z} \cup \{p\}^{k+1}$
with that of~${\cal Z} \cup \{p\}^{k}$ from explicit identification of the expected linear obstruction
schemes is presented.  
Particularly, after intersecting such a scheme with a hyperplane~$H$ we produce a collection~${\cal W}$
of multiple points of~$H$ so that
$$h_{\proj^n}({\cal Z} \cup \{p\}^{k+1},m) - h_{\proj^n}({\cal Z} \cup \{p\}^k,m) \leq
{{n+k-1} \choose {n-1}} - h_{H}({\cal W},k).$$
Further, equality occurs exactly when the scheme~${\cal Z} \cup \{p\}^{k+1}$ has only the expected linear 
obstructions given by~${\cal Z} \cup \{p\}^k$ in the relevant degree~$m$, as in Definition~\ref{exp}.

From Lemma~\ref{key} we obtain the Main Theorem (Section~6) on describing the Hilbert function of a collection of multiple
points.  Namely,  an upper bound on the Hilbert function of a collection of multiple
points of~$\proj^n$ is obtained from evaluation of the Hilbert function of collections of 
fat points in~$\proj^{n-1}$, as follows:

\begin{theorem} \label{ubda} Let~$n, m, d \in \nats$ and~$A=(k_1, \ldots, k_d) \in \nats^d$.
Take~$C_{ji}=(k_1, \ldots, k_j)+\overline{i-m}$ for each~$j=1, \ldots, d-1$, and~$i=0, \ldots, k_j-1$.

For each~$A$-subscheme~${\cal Z} \subset \proj^n$ there are naturally induced~$C_{ji}$-subschemes of~$\proj^{n-1}$,
$W_{ji}$, so that
$$h_{\proj^n}({\cal Z}, m) \leq \deg {\cal Z} - \sum_{j=1}^d \sum_{i=0}^{k_j-1} h_{\proj^n}(W_{j,i},i).$$
Equality holds if and only if only the expected linear obstructions to~$Z$ occur in degree~$m$.
\end{theorem}

This result has  immediate implications toward comparing the Hilbert function of a  fat point scheme
with the 
function~$G(d,A,n)_m$ proposed by Iarrobino, as seen by the basic  properties of the function~$G$.
From this we obtain the conclusion of Theorem~\ref{siwa}.

\medskip
 
In Section~8 we examine general applications of Theorem~\ref{ubda} toward~$\proj^n$.
We restrict attention mainly to cases in which we are ``one step away'' from the expectation
of maximal rank: only (multiple) lines are predicted to appear
as the positive-dimensional schemes  in the base locus of the given linear systems, in the
formulation of the conjectural function.
We obtain in Corollary~8.8, for example, a geometric analogue of an algebraic result of Iarrobino (Section~4).

As a further application of the methods developed here,
 we verify the Strong Conjecture in two main settings (Section~9).
The first of these involves an hypothesis that the sum of 
multiplicities is not too large compared to the degree~$m$ and
the dimension~$n$ (which, technically, may be viewed as asserting
that a rational normal curve cannot possibly impede the Hilbert
function).  This is obtained from Theorem~\ref{ubda} (for the
upper bound) together with Castelnuovo techniques (lower bound).
[An experienced reader might be shocked to find the first 
appearance of said Italian mathematician occurring so late in
a paper of this author!]

It follows then that we may narrow down the possible exceptions to the
Strong Fr\"oberg-Iarrobino hypothesis.  
For example (as we generalise in Corollary~\ref{n++3}):

\begin{cor}\label{nplus3}
Let $n,m \in \nats$.  
Let~$A = (k_1, \ldots, k_{n+3}) \in \nats^{n+3}$.
Then a generic~$A$-scheme ${\cal Z}$ does satisfy the Strong  Fr\"oberg-Iarrobino Hypothesis
in degree~$m$ provided that
$$ \sum_{i=1}^{n+3} k_i \leq mn+1; $$
namely, for the rational normal curve~$C$ through the~$n+3$ points,
the degree of~$C \cap {\cal Z}$ is  at most the value of the
Hilbert function of the curve in this degree.
\end{cor}

The second case studies instances of points of equal multiplicity~$k$
with focus on the case of degree~$k+1$ (the first nontrivial case).
We showed in~\cite{me} (using, of course, Castelnuovo methods) that the
value proposed by Iarrobino gives a {\it lower bound} in certain cases
(much to the consternation of Iarrobino at the time, since this is the
``difficult part''!).  Here we simply apply Theorem~\ref{ubda} to find equality.

However, in Section~10 we find counterexamples to the Strong Conjecture in~$\proj^n$ 
for each~$n=4, 5, 6$ based on intersection with a rational normal curve.
We comment, then, on determining when the Strong Hypothesis should apply to
fat point schemes (homogeneous or otherwise).

Finally, we remark on the algebraic conjectures of Fr\"oberg, equipped with
the extra information from the geometric viewpoint.
We describe how Conjecture~\ref{conj1} implies the ``Koszulness'' of the minimal resolution of an ideal 
generated by powers of linear forms, and whence, general forms.

In sum, the geometric evidence presented here gives structure to %
 the Fr\"oberg-Iarrobino conjectures.
  Indeed, the Weak Conjecture
appears tractable technically.  Furthermore one needn't regard 
the Weak Conjecture as the ``second best'' result to obtain
along these lines, but as a first step toward verifying
the strong conjecture (and evaluating exceptional cases).
Namely, from the Weak Conjecture it would follow that equality
of the Hilbert function of an~$A$-subscheme of~$\proj^n$  and~$G(d,A,n+1)$ is an open condition (on $(\proj^n)^d$),
and hence
may be verified by {\it producing} a scheme exhibiting such equality.
Particularly, a usual strategy for verifying upper bounds on the Hilbert function
of a general scheme is to construct a scheme that satisfies these conditions.
In the situation of maximal rank, this always suffices; but it is necessary here to have a lower bound.
Moreover, as we see in Theorem~\ref{siwa}, the Weak Conjecture allows us to characterise schemes that
do obey the Strong Hypothesis (including the Strong Conjecture). 


The structure of the paper is as follows.
We fix notation in Section~2.
Next we consider the context of the problem at hand.
We begin in Section~3 by describing basic results and techniques that
may be used to  predict maximal rank.
Then in Section~4 we present the conjectures of Fr\"oberg 
and Iarrobino, which imply in particular how maximal rank cannot
always be achieved.
From basic observations on intersection multiplicity we obtain in Section~5 the 
geometric interpretation of these conjectures along with refinements of the conjectures.  

In Section~6 we prove the main theorem and present its connection
to the conjectures. 
Then, in Section~7, we validate the 
 Weak Conjecture in~$\proj^3$.
Further consequences in terms of verifying cases of the Weak Conjecture in~$\proj^n$
are given in Section~8.

In Section~9 we  verify some instances of the Strong Conjecture
{\it using} the acquired information on the weak one.
We exhibit in Section~10 counterexamples to the Strong Conjecture in~$\proj^n$ for~$n \leq 6$,
and remark on revision of this conjecture and

We end in Section~11 with discussion and speculation on the 
algebraic versions of the strong conjecture.

\medskip

\noindent{\bf Acknowledgements.}  I am very grateful to Tony Iarrobino
for his inspiring my work along this direction and for discussions.
Thanks also to Joe Harris for sharing his comments.

\section{Basic Bookkeeping}

Given a projective subscheme~${\cal Z} \subset \proj^n$ we write
${\cal I}_{{\cal Z}, \proj^n}$ (or ${\cal I}_{\cal Z}$)
 for its ideal sheaf.

Given a reduced subvariety~$X \subset \proj^n$  and~$a \in \ints$,
 define~$X^a$ as the
subscheme~${\cal Z} \subset \proj^n$ defined by~${\cal I}_X^a$.
(So in  case~$a \leq 0$ we have~$X^a = \emptyset$.)

Since we shall consider collections of fat points of various 
multiplicities, let us keep some of the bookkeeping straight as follows:

Given~$A=(a_1, \ldots, a_d) \in \ints^d$ 
define an~$A$-{\bf subscheme of}~$\proj^n$
as a union~${\cal Z}$ of~$\{p_1\}^{a_1} \cup  \ldots \cup \{ p_d\}^{a_d}$
where~$\{p_1, \ldots, p_d\}$ is a set of~$d$ points in~$\proj^n$.
In case~$A=(a, \ldots, a)$ (and the number of points~$d$ is made clear)
we shall write~$A=\overline{a}$ and refer to an~$A$-scheme as a
{\bf homogeneous} subscheme.
We say that~${\cal Z}$ is a {\bf generic $A$-subscheme}  (or {\bf general}) if
the set~$\{p_1, \ldots, p_d\}$ is generic.

As in~\cite{iar} we define:

\begin{defn}
Denote by $HPTS(d, (k_1, \ldots, k_d),n+1)_m$
 the value of Hilbert function in
degree $m$ of the generic union of $d$ points in $\proj^n$ of
multiplicities $k_1, \ldots, k_d$.
\end{defn}
Let us extend this definition in the obvious manner
to the situation of perhaps having negative entries in the uple~$A$:
$$HPTS(d, A, n+1)_m := HPTS(d, (\max(a_1, 0), \ldots, \max(a_d,0),n+1).$$

For such~$A$ we write~$|A|=\sum_{i=1}^d \max(a_i, 0)$ 
and~${\ell}(A)= \#  \{1 \leq i \leq d: a_i > 0\}$.

Given~$A, B$ both uples, we say that~$A$ and~$B$ are equivalent, or
that~$A$ {\bf may be written as~$B$} if an~$A$-scheme is a~$B$-scheme.
This says that~$\ell(A)=\ell(B)$ and if~$a_1, \ldots, a_d$
and~$b_1, \ldots, b_d$ are the positive entries of~$A$ and~$B$
respectively ($d=\ell(A)$) then there exists a permutation~$\sigma$
on~$d$ letters
so that~$a_i = b_{\sigma(i)}$ for~$i=1, \ldots, d$.

Now take~$A=(a_1, \ldots, a_d) \in \nats^d$ so that~$a_i > 0$
for~$i=1, \ldots, d$.

Given~$B \in \nats^r$ we shall say that~$B \subseteq A$ 
if~$B=(a_{i_1}, \ldots, a_{i_r})$ with~$1 \leq i_1 < \ldots < i_r \leq d$
(so~$B$ does respect ordering).

If~$B \in \ints^r$ we say that~$B \leq A$ if~$B$ may be written 
as~$(b_1, \ldots, b_d) \in \nats^d$ with each~$b_i \leq a_i$, and that
$B < A$ if~$B \leq A$ but~$B$ may not be written as~$A$.

\section{General Context}
We  pay particular attention throughout this paper to the Hilbert function of a generic collection
of multiple points in projective space as a start toward the study of arbitrary
collections of fat points.
In each degree~$m$ the Hilbert function~$HPTS(d,A,n+1)_m$ of a (general)
$A$-scheme is bounded above by the degree of the scheme, with equality for
$m$ sufficiently large.
Let us describe here some of the circumstances under which the two 
quantities are known (or conjectured) to agree.
We refer the reader as well to the very readable and detailed accounts in the surveys of Ciliberto~\cite{cil}
and Harbourne~\cite{h2}, particularly for thorough descriptions of work on~$\proj^2$.

Let~$A=(k_1, \ldots, k_d) \in \nats^d$ and
consider the Hilbert function of an~$A$-subscheme of~$\proj^n$.
For convenience, let us assume that~$k_1 \geq k_2 \geq \ldots \geq k_d$, for now.

When $d=1$ we have
$$HPTS(1,\bar{k},n+1)_m = 
\min \left( {{n+m} \choose m}, {{n+k-1} \choose {k-1}}\right);$$
that is, the scheme has the maximal rank property.

Now consider~$d \geq 2$.  For each~$m \leq k_1+k_2-2$ 
we have
$$HPTS(2,(k_1,k_2),n+1)_m < {{n+k_1-1} \choose 2}+{{n+k_2+1} \choose 2}. $$

Indeed, if~$ m \leq k_1+k_2-2$, look at a~$(k_1, k_2)$-scheme~$\cal Z$
and take the line~$L$ spanned by the two reduced points on~$\cal Z$.
Then~$\deg L \cap Z = k_1+k_2$.  However, $L$ imposes only~$m+1$ conditions
on the linear system of~$m$-ics.  Hence~${\cal Z} \cap L$ cannot
impose independent conditions on~$m$-ics, so neither does~$\cal Z$.
We regard this as an~{\bf expected linear obstruction}.
(We shall make this notion precise in Definition~\ref{exp}.)

Likewise, extending to~$A=(k_1, \ldots, k_d)$ we have that an~$A$-scheme
cannot have maximal rank unless
$HPTS(d,A,n+1)_{m} = {{n+m} \choose n}$
for~$m=k_1+k_2-2$
(that is, the degree of an~$A$-scheme is large enough).

In the case of degree~$m \geq k_1+k_2-1$, 
the m\'ethode d'Horace of ~\cite{h} does (essentially) 
apply toward verifying inductively a given  case of maximal rank
from ones occurring in  lower degree and in lower dimension.
This led to the following 
asymptotic result of J. Alexander and A. Hirschowitz:

\begin{theorem}(\cite{ah5})
Given~$n,k \in \nats$, there is a quantity~$d(n,k)$ so that for all~$d \geq d(n,k)$,
$$HPTS(d,\overline{k},n+1)_m = \min( d {{n+k-1} \choose n}, {{n+m} \choose m}).$$
\end{theorem}

One would like to sharpen this to an actual upper bound for~$d(n,k)$;
indeed, one
that should be independent of the
multiplicity~$k$.
Equivalently, the theorem predicts that there is a value~$c(n,k)$
so that whenever~$m \geq c(n,k)$ we have equality between~$HPTS$
and  the desired quantity.
Again one should like a bound on such a value, independent of~$n$.
This asks for results in cases~$m \leq 2k-2$ (in which maximal rank
cannot be achieved) to 
obtain a starting point in the induction process.
Unfortunately,  
the m\'ethode diff\'erentielle does not apply well to the
low degree cases~($m \leq 2k-2$).  

On the bright side (as used, e.g.,  in~\cite{a}, \cite{me2}, \cite{methree}) to proceed by induction
on degree one does not require that {\it every } collection of
multiple points has maximal rank.  
For example,  to show that
in  degree~$m \geq 3k$ a generic collection of~$k$th order points 
exhibits maximal rank in degree~$m$, it suffices to concoct 
a fat point scheme of large enough degree imposing independent conditions in
degree~$m-k+1$.  But the degree required is strictly (and significantly)
less than than the boundary value of~$\displaystyle {{n+m-k+1} \choose n}$.
With any luck, then, down-to-earth methods in low degrees (e.g. \cite{me})
should yield $k$-schemes of maximal rank in degree~$2k-1$, say, and
of large enough degree to proceed inductively.

\begin{conj}
Let~${\cal Z} \subset \proj^n$ be a generic collection of points
of multiplicities at most~$k$.  Then for each~$m$ with~$m \geq 3k$
we have
$$h_{\proj^n}({\cal Z}, m)= \min(\deg {\cal Z}, {{n+m} \choose m}).$$
\end{conj}

Rewriting in terms of the number of points, we make (or, at least, estimate) the following:


\begin{conj}
A generic collection of~$d$ fat points of~$\proj^n$ has
maximal rank if~$d > 2^{n+1}$.
\end{conj}

The case of the projective plane has seen a good deal of progress in recent years (see \cite{h2}).

First, Nagata's method \cite{n}  of blowing up multiple points led
to his observation that the ``visible'' obstructions to maximal rank
occur precisely when there are at most eight multiple points.
Hence the main information on Hilbert functions of fat points in~$\proj^2$
should be gleaned from the study of linear systems on Hirzebruch surface.
Under this insight, he evaluated (for example):

\begin{theorem} \cite{n} Let~$d, m \in \nats$, with $d \leq 6$.  
Let~$A=(k_1, \ldots, k_d) \in \nats^d$ so that~$k_1 \geq k_2 \geq \ldots \geq k_d$
and~$ \sum_{i=1}^5 k_i \leq 2m+1$. 
Then a generic~$A$-subscheme of~$\proj^2$ has maximal rank in degree~$m$.
\end{theorem}

Notice that the assumption on~$\sum_{i=1}^5 k_i \leq 2m+1$ ``prevents'' the
conic through the first five points from interfering with maximal rank.
Quite generally, efforts of Segre 
on this theme yield the following:

\begin{conj}\label{hh} (\cite{se})  If a collection~${\cal Z}$ of multiple points
 in $\proj^2$ does not display maximal rank in degree~$m$ then there
is a curve~$C$ in the base locus of the system of~$m$-ics through~${\cal Z}$
for which
$$\deg {\cal Z} \cap C > h_{\proj^2}(C,m).$$
\end{conj}

Harbourne and Hirschowitz (\cite{h}, \cite{h1}) conjectured on the Hilbert function of a collection
of multiple points as well,
referring to the induced linear systems on the blow-up of~$\proj^2$ with respect to the points.
Hirschowitz used his  ``m\'ethode d'Horace'' to evaluate the Hilbert
function of a generic union of double points and likewise for triple points of~$\proj^2$.
Harbourne has made considerable refinements on the  predictions of Nagata (e.g., \cite{h1}) on determining exceptional
cases to maximal rank, with support on eight or fewer points; moreover, attacking the more difficult 
problem of finding free resolutions of the ideal of such schemes.  
He has also developed an algorithm for computation of
the Hilbert function of fat points in~$\proj^2$ 
(see \cite{h3}, along with the computer programme running on his web page).

Ciliberto and Miranda \cite{cm1}, \cite{cm2} have made further
significant progress toward this problem.
They show (\cite{cms}) that Segre's conjecture indeed implies those of Harbourne and Hirschowitz. 
Their method involves passing between~$\proj^2$ and the Hirzebruch surface~$\f_1$.
This led to a complete analysis of the cases~$HPTS(d, (k_1, \ldots,k_d),3)$, $k_i \leq 3$
and~$HPTS(d,\bar k,3)$, $k \leq 12$. 
They, along with Orecchia, applied the technique later to extend to fat points of
(equal) multiplicity up to~20 \cite{cmo}.
Further, Ciliberto and Miranda, proceeding on the description in the conjectures of~\cite{h} and~\cite{h1}
 give a precise identification of how curves may inhibit the Hilbert function of a collection
of multiple points in the plane, yielding nice refinements in~\cite{cm1} 
 of the Segre conjecture.

Meanwhile, up in higher dimensions we do have the results of Alexander and Hirschowitz
(\cite{h}, \cite{a}, \cite{ah1}, \cite{ah2}, \cite{ah3}.) on the Hilbert function of generic double points in projective space
using (nontrivial!) variations on the m\'ethode d'Horace.  
(To complete that story, the author evaluated the one missing case \cite{me4}!)

Next, the author shows in~\cite{methree} that a generic union of double and
triple points in~$\proj^n$ does exhibit maximal rank in degree at least~7, 
and exhibits schemes that do not satisfy the Strong Hypothesis in lower degrees.
This is done by expanding on the simplified version of ``Horace'' given in~\cite{me2}
for a ``brief proof'' of the Alexander-Hirschowitz result.
In the latter paper, this appears as an organisational means toward the result
at issue, whereas in multiplicity~three its use is critical, not only in bookkeeping 
but in bypassing assumptions on, say, Harbourne-Hirschowitz-type conjectures.

\section{Fr\"oberg-Iarrobino conjectures}

The conjectures of Fr\"oberg and Iarrobino apply to situations
in which the Hilbert function of a collection of multiple points
{\em should not} attain maximal rank in a given degree; namely,
when that degree is small compared to to multiplicities.

We present here a background on these conjectures, together
with pertinent known results.  First, the conjectures of Fr\"oberg on
the behaviour of an ideal generated by a generic set of forms in a
(graded) polynomial ring are described.  Next, we look at the
``specialisations'' made by Iarrobino to an ideal generated by 
powers of linear forms, chosen generically. We see then the
interpretation of Emsalem and Iarrobino of the latter problem to
the Hilbert function of multiple points.

\subsection{ Fr\"oberg conjectures on the ideal of general forms}$\mbox{~}$ \\

Let us take $R=K[Y_0, \ldots, Y_n]$, as a graded ring.

Fr\"oberg \cite{f} considers an ideal~$J$ generated by a generic collection
of forms of specified degrees.  His conjectures suggest that the
minimal free resolution of~$J$ in ``looks like'' that
of a complete intersection until the the degree is, by numerics, 
expected to be large enough that the ideal should contain all 
forms of that degree.

\begin{defn}
Let $A=(j_1, \ldots, j_d) \subset \nats^d$.
Let $J$ be an ideal generated by a  general set of~$d$ forms
in~$n+1$ variables, of degrees~$j_1, \ldots, j_d$.
 
We write
$HGEN(d,A,n+1)_m=\dim R_m - \dim J_m$.
\end{defn}

\noindent{\bf Example.} $d \leq n+1$.  Then $J$ is a complete
intersection ideal with minimal free resolution:

$$0 \rightarrow R(-j_1-\ldots-j_d) \rightarrow \ldots \rightarrow
\oplus_{1 \leq i < k \leq d} R(-j_i-j_k) \rightarrow 
\oplus_{i=1}^d R(-j_i) \rightarrow J \rightarrow 0.$$
Hence
$$\dim R_m - \dim J_m = \sum_{B \subseteq A} 
(-1)^{\ell(B)} {{n+m-|B|} \choose n}.$$
For example, in the case $j=j_1=\ldots=j_d$ we have 
$$\dim R_m - \dim J_m = 
\sum_{t=0}^n (-1)^{t} {d \choose t} {{n+m-tj} \choose n}.$$

Then, for $d \geq n+2$, Fr\"oberg conjectures that, taking $J$ generated
by general forms, then $\codim J_m$ should act as if the only
relations between generators are Koszul for as long as possible, namely,
until the first (numerical) opportunity for $J_m = R_m$.
Specifically:

\begin{defn}  Define a function~$F$ as follows:  
Given~$n, d \in \nats$, $A \in \nats^d $, 
and~$m \in \nats$, set
$$F'(d,A,n+1)_m = \sum_{B \subseteq A} (-1)^{\ell(B)} {{n+m-|B|)} \choose n},$$ 
and 
$$ F(d,A,n+1)_m = 
\begin{cases}
 0, & \mbox{~if~}
F'(d',A',n+1)_m  \geq {{n+m} \choose m} \\
& {~for~some~} d' \leq d ~\mbox{~and~} A' \subseteq A \\
F'(d,A,n+1)_m, & \mbox{otherwise}. 
\end{cases}$$
\end{defn}

Then the conjectures of Fr\"oberg are:
 
\begin{conj}
Strong Fr\"oberg Conjecture(SF): $$HGEN(d,A,n+1)_m=F(d,A,n+1)_m.$$
\end{conj}
 
\begin{conj}
Weak Fr\"oberg Conjecture(WF): $$HGEN(d,A,n+1)_m \geq F(d,A,n+1)_m.$$
\end{conj}

Under the following hypotheses the strong conjecture is known to hold 
(see \cite{iar} for more details): $d \leq n+1$ (as we have just
observed), $d=n+2$ ( R. Stanley, \cite{s}), 
$n=1$ (R. Fr\"oberg, \cite{f}), and~$n=2$ (D. Anick, \cite{an}).

In the case of equal degree, that is, $A=(j, \ldots, j)$, the
strong conjecture has been verified in the following cases:
$m=j+1$ (M. Hochster and D. Laksov, \cite{hl}), 
$n \leq 10$ and $j \leq 2$, along with $n \leq 7$ and $j \leq 3$
(Fr\"oberg and J. Hollman, \cite{fh}); and the more detailed 
hypotheses of M. Aubry (\cite{au}), that~$m$ is sufficiently close to~$j$.

Iarrobino shows that 
in the ``first Koszul interval'' for a specified case, the
Weak Fr\"oberg Conjecture  may be verified by a lower degree case of
the Strong Fr\"oberg Conjecture:

\begin{theorem}(\cite{iar})\label{ith}
Take~$A=(j, j_1, \ldots, j_d)$ and~$C=(j_1, \ldots, j_{d})$.
Assume that $j \leq \min \{j_i\}$ and $2j \leq m \leq 3j$. 

If $$HGEN(d,C,n+1)_{m-j}=F(d,C,n+1)_{m-j}$$ then
$$HGEN(d,A,n+1)_{m} \geq F(d,A,n+1)_{m}.$$
\end{theorem}

We shall obtain in Corollary~\ref{pith} a geometric analogue of this result via Theorem~\ref{ubda}.

\subsection{ Iarrobino's conjectures on general linear forms} {$\mbox{~}$ }\\

Iarrobino extends the conjecture to powers of linear forms in
appropriate characteristic:

\begin{defn}
Let $A=(a_1, \ldots, a_d)$.
Let $J$ be an ideal generated by a  generic set of~$d$ powers of
linear forms
in~$n+1$ variables, of degrees~$a_1, \ldots, a_d$.

Write
$HPOWLIN(d,A,n+1)_m=\dim R_m - \dim J_m$.
\end{defn}

Of course, by upper-semicontinuity, 
$$HPOWLIN(d,A,n+1)_m \geq HGEN(d,A,n+1)_m.$$ 

This yields:

\begin{conj}\label{slin} Strong Algebraic Fr\"oberg-Iarrobino:
Let~$n,m,d, a \in \nats$.
Then
$$HPOWLIN(d,\bar{a},n+1)_m \geq F(d,\bar{a},n+1)_m,$$ 
Further, let $p$ be the characteristic of $ \mathcal K$.
If $p=0$ or $ m > p$, we have
$$HPOWLIN(d,\bar{a},n+1)_m=F(d,\bar{a},n+1)_m,$$ 
except in the following circumstances:
$d =n+3$, $d=n+4$,  $n=2$ and $d = 7$ or $8$; $n=3$ and $d=9$.
\end{conj}

\begin{conj}
Weak Algebraic Fr\"oberg-Iarrobino:
$$HPOWLIN(d,A,n+1)_m \geq F(d,A,n+1)_m.$$
\end{conj}


{\bf Remark:} We shall see more explicitly how the exceptional values arise from the
interpretation to the postulation of multiple points described next.

\subsection{ Conjectures on multiple points} \mbox{~}\\

Emsalem and Iarrobino deduced
from Macaulay duality that $HPOWLIN$ is related to the
Hilbert function of multiple points
in appropriate characteristic.
Namely,
\begin{theorem}(\cite{ei})  Let $A=(k_1, \ldots, k_d)$
 and $A'=(m+1-k_1, \ldots, m+1-k_d)$.
If~$\cchar {\cal K} = 0$ or
$\cchar {\cal K} > \max(m,k_1, \ldots, k_d)$ then
$$HPTS(d,A,n+1)_m = \dim R_m - HPOWLIN(d,A',n+1)_m.$$
\end{theorem}

To see the main idea of this theorem, 
let us look at the case~$\cchar {\cal K}=0$.
Let~$S={\cal K}[X_0, \ldots, X_n]$ 
and~$R={\cal K}[{\partial} / {\partial} X_0, \ldots, 
{\partial}/ {\partial }X_n]$ (with~$R$ regarded as a polynomial ring
in the ``dummy variables'' ${\partial/\partial X_i}$). 

Macaulay duality refers to the obvious perfect pairing
${\Phi_m} : R_m \times S_m \rightarrow {\cal K}$. 
(Of course, it is not quite natural in that it depends 
on the co\"ordinate choice.)
Given an ideal~$I$ of~$S$ we obtain, in each degree~$m$,
$$I_m^{\perp}=\pi_{R}(\ker \Phi |_{R_m \times I_m}) \subset R_m,$$
so that~$\dim I_m^{\perp} = \dim S_m - \dim I_m$.
Certainly for two ideals~$I, J$ we have
$(I \cap J)_m^{\perp} = I_m^{\perp} + J_m^{\perp}$.

The case in point given by Emsalem and Iarrobino  starts with
the ideal of a multiple point, say~$I=(X_1, \ldots, X_n)^k \subset S$.
Then~$I_m = (X_1, \ldots, X_n)^k S_{m-k} $ 
so that~$I_m^{\perp}={\partial/\partial X_0} R_{k-1}$.
(Each monomial, say, $M$~in~$I_m$ 
satisfies~${\partial/\partial X_0}^{m-k+1}M=0$.)

Whence the dual to the~$m$th graded piece of the ideal 
of a~$(k_1, \ldots, k_d)$-scheme is the $m$th piece of an ideal 
generated by~$d$ powers of linear forms, 
the powers being~$m+1-k_1, \ldots, m+1-k_d$.

So, Iarrobino's conjectures on powers of linear forms ``translate'' to
 candidates for values of the function $HPTS$ evaluating the Hilbert function
of fat points.  Moreover, note that in the latter setting the characteristic
of the field need not play a r\^ole.
This yields the following:

\begin{defn}\label{G}
Let $$G'(d,(k_1, \ldots, k_d),n+1)_m = \dim R_m -F'(d,(m+1-k_1, \ldots, m+1-k_d), n+1)_m,$$
and
$$ G(d,(k_1, \ldots,k_d),n+1)_m = 
\dim R_m - F(d, (m+1-k_1, \ldots, m+1-k_d), n+1)_m.$$
\end{defn}

The geometric versions of the Fr\"oberg-Iarrobino conjectures  are then
the (stronger) conjectures:

\begin{conj}\label{wpt}
Weak Fr\"oberg-Iarrobino (WFI): $HPTS(d,A,n+1)_m \leq G(d,A,n+1)_m$.
\end{conj}

\begin{conj}\label{spt}
 
Strong Fr\"oberg-Iarrobino (SFI: homogeneous): 
For each~$n,d,m \in \nats$ we have $HPTS(d, \overline{k}, n+1)_m \leq G(d,\overline{k}, n+1)_m$.
Further, we have $HPTS(d,\overline{k},n+1)_m=G(d,\overline{k},n+1)_m$,
except perhaps when one of the following conditions holds: 
$d =n+3$, $d=n+4$,  $n=2$ and $d = 7$ or $8$; $n=3$, $d=9$, $m=2k$; or
$n=4$, $d=14$, $m=2k$ and $k=2$ or $3$.
\end{conj}

Hence in each {\it homogeneous} case $d \geq n+5$ we have
$$ SFI \implies SF \implies WF \implies WFI.$$

As previously stated, we shall say a generic~$A$-subscheme of~$\proj^n$ (possibly of mixed multiplicities)
 satisfies the Strong Fr\"oberg-Iarrobino 
Hypothesis if its Hilbert function agrees with the conjectural value~$G(d,A,n+1)$.
(Similarly, we  say that it satisfies the Strong Hypothesis in a given degree if the value of
the Hilbert function in this degree is equal to the specified value. 
Likewise, we refer to an arbitrary scheme as satisfying the Strong Hypothesis if its Hilbert function
agrees with the value of the corresponding function~$G$.)

It is straightforward to verify the Strong Conjecture in the case of~$d \leq n+1$
fat points in~$\proj^n$ by examining the intersection ideal.
A geometric proof appears in Section~9.
For~$d=n+2$ the conjecture holds valid as well, again by~\cite{s}.

On the other hand, the exclusion of cases $d =n+3, n+4$, 
and so forth in Conjecture~\ref{slin} and then Conjecture~\ref{spt} do arise from
the expectation of special varieties through points that inhibit
the Hilbert function. 

We began a direct geometric analysis of the conjectures on homogeneous $\overline{k}$-subschemes of~$\proj^n$ in \cite{me}.
Since the cases of degree~$m \leq k$ are trivial to verify (at most~$n+1$ points are involved) the focus there
is finding criteria determining the inequality: $HPTS(d,k,n+1)_{k+1} \geq G(d,k,n+1)_{k+1}$.
The result (see Proposition~\ref{plus1}) invokes the identification of neighbourhoods of planes of each dimension lying in
in the base locus of a $(k+1)$-ic linear system through a~$\overline{k}$-scheme (made precise in Section~5 below).

Here we expand on this method for use in the general setting.

\section{Interpreting conjectures geometrically}

We describe the correspondence of the Fr\"oberg-Iarrobino conjectures
with the issue of linear obstructions to the maximal rank of a collection
of fat points.
In particular, we make precise the notion of a scheme that displays
only the expected linear obstructions (as in Theorem~\ref{ubda}),
which we ``expect'' to satisfy  the Strong Hypothesis.
We start with a simple identification of multiple planes (of each dimension)
 that must appear in the
base locus of a linear system through a given collection of multiple points.
We compare this information
to the conjectures at hand by means of intersection degrees.
This leads naturally to extensions and refinements of these conjectures. 

{\bf Basic Observation.}
Let $p, q \in \proj^n$ and $L=\spax\{p,q\}$.
By B\'ezout (or by Lagrange-Hermite!), any $m$-ic form 
vanishing on $\{p\}^k \cup \{q\}^j$
vanishes on~$L$ if $m\leq k+j-1$.
Further, such an~$m$-ic must vanish on~$L^{k+j-m}$.
This is easy to see in characteristic~0, simply by taking derivatives.
Generally, take $I=(X_0, X_1, \ldots, X_n)^k \cap (X_0, X_2 \ldots, X_n)^j$.
Suppose that~$F \in I_m$ and write~$F=X_0^r X_1^s G$ so that neither~$X_0$ nor~$X_1$
divides~$G$.  Let us write~$G=N+G_1$ so that~$N$ is a monomial 
(divisible by neither~$X_0$ nor~$X_1$) and~$G_1$ is a sum of fewer monomial terms than~$G$.
We have that  $X_1^s G \in (X_1, \ldots, X_n)^j$, so that $\deg G \geq j-s$, and
hence $r \leq m-j$.  Likewise, $s \leq m-k$, so~$\deg G \leq k+j-m$.
In particular~$N \in (X_2, \ldots, X_n)^{k+j-m}$, so that we may replace~$G$ by~$G_1$
and repeat the procedure, reducing the number of terms at each step.
Whence~$G$ vanishes on~$L^{k+j-m}$ and so does~$F$.

\medskip

Consequently, we obtain:

\begin{lemma}\label{lag}
Let~$P, Q \subset \proj^n$ be projective subspaces.  
Then every~$m$-ic vanishing on~$P^k \cup Q^j$ must vanish on
$\spax(P \cup Q)^{k+j-m}$.
\end{lemma}

{\it Proof:} Suppose that~$F \in I(P^k \cup Q^j)$.  
Let~$t \in \spax(P \cup Q)$.   Then~$t \in L$ where~$L$ is the line
between two points, say, $p, q \in P \cup Q$. 
As we have just seen, $F \in I(t)^{k+j-m}$
since~$F \in I(p)^k \cap I(q)^j$.
\hfill $\Box$.

\medskip

This yields the following generalisation of our initial note:

\begin{cor}\label{lagp}  Let $n,m, d \in \nats$ and $(k_1, \ldots, k_d) \in \nats^d$.
Suppose that~$P_1, \ldots, P_d$ are projective subspaces of $ \proj^n$.
Then every~$m$-ic through $P_1^{k_1} \cup \ldots \cup P_d^{k_d}$
must vanish on $\spax( P_1 \cup  \ldots \cup  P_d)^r$ for 
$r = k_1 + \ldots + k_d - (d-1)m$.
\end{cor}

{\it Proof:} Apply induction together with Lemma~\ref{lag} to the union of 
$P_d^{k_d}$ and~$\spax (P_1 \cup \ldots \cup P_{d-1})^j$,
where~$j = k_1+\ldots+k_{d-1} - (d-2)m$.
\hfill $\Box$.

\medskip

This motivates the following description: 

\begin{defn} \label{exp}  Let $n,m,d, k \in \nats$, and~$(k_1, \ldots, k_d) \in \nats^d$.

Let~${\cal Z} = P_1^{k_1} \cup \ldots \cup P_d^{k_d}$ and~$P \subset \proj^n$,
where~$P_1, \ldots, P_d,$ and~$P$ are projective subspaces of~$\proj^n$
(such as points). 

Let us refer to the {\bf expected linear obstruction scheme on $P^k$} induced by
$m$-ics through $P^{k-1} \cup {\cal Z}$ as the subscheme of~$P^k$ predicted
by Corollary~\ref{lagp} namely,  
$$(P^{k-1} \cup \bigcup 
\spax(P \cup  P_{i_1} \cup \ldots \cup P_{i_r})^{j(i_1, \ldots, i_r)}) \cap P^k,$$ 
where the union ranges over sets of indices~$(i_1, \ldots, i_r)$ 
for which~$1 \leq i_1 < i_2 < \ldots < i_r \leq d$, and
$j(i_1, \ldots, i_r)= k+ k_{i_1}+ \ldots+ k_{i_r}-rm-1$. 

We shall say that~$P^k \cup \z$ is 
{\bf only linearly obstructed} by~$P^{k-1} \cup \z$ in degree~$m$  
if the base locus of the linear system of~$m$-ics
 through~$P^{k-1} \cup \z$ meets~$P^k$ in precisely the expected linear obstruction scheme. 
 
We say that~$P^k \cup \z$ is {\bf only linearly obstructed} by $\z$  
in degree~$m$ if $P^{\ell} \cup \z$ is only linearly obstructed 
by~$P^{\ell-1} \cup \z$  for each~$\ell$ with~$1 \leq \ell \leq k$. 
 
Finally, we say~$\z$ has {\bf only the expected linear obstructions} in degree~$m$ 
if for each~$1 \leq j \leq d$ we have 
that~$P_1^{k_1} \cup \ldots \cup P_j^{k_j} \cup \ldots \cup P_d^{k_d}$  
is only linearly obstructed 
by~$P_1^{k_1} \cup \ldots \cup P_j^{0} \cup \ldots \cup P_d^{k_d}$.

\end{defn}
 
Note that, in the context of the definition, if the set of planes 
is generic, the planes are of equal dimension,
 and~$k_1 \geq k_2 \geq \ldots \geq k_d$, 
then~$\z$ is nonlinearly obstructed if and only 
if~$\z$ is nonlinearly obstructed by $P_2^{k_2} \cup \ldots \cup P_d^{k_d}$.  
 
\medskip

Let us now make an obvious and (hence!) useful simplification to Definition~\ref{exp}.

\begin{lemma}\label{lobs} Let $n,m,k \in \nats$ and $(k_1, \ldots, k_d) \in \nats^d$.
Take projective subspaces~$P, P_1, \ldots, P_d$ of $\proj^n$ and let
${\cal Z}=P_1^{k_1} \cup \ldots \cup P_{d}^{k_d}$.

Take $\Upsilon$ as the expected linear obstruction subscheme of~$P^k$ given
by~$m$-ics through~$P^{k-1} \cup {\cal Z}$ and $\Upsilon_1$ as the 
scheme~$\Upsilon_1=(P^{k-1} \cup \bigcup_{i=1}^d \spax(P \cup P_i)^{k+k_i-m-1}) \cap P^k$.
Then~$\Upsilon=\Upsilon_1$.
\end{lemma}

{\it Proof:} We do have~$\Upsilon_1 \subseteq \Upsilon$.

After intersection with a hyperplane we reduce to having~$P=\{p\}$, $p \in \proj^n$.

By symmetry, it suffices to see that for each~$r \leq d$ and~$j=k+k_1+\ldots+k_r-rm-1$
we have~$(\{p\}^{k-1} \cup \spax(p \cup P_1 \cup \ldots \cup P_r)^j) \cap \{p\}^k \subseteq \Upsilon_1$
.
Let us fix such an $r$, then.

Choose an open affine subspace of $\proj^n$ containing~$p$ (viewed as an origin point)
and take~${\frak m}$ as the maximal ideal of~$p$.
In the associated projective space the quotient~${\frak m}^{k-1}/{\frak m}^{k}$ identifies
forms of degree~$k-1$.  Then the desired inclusion follows straight from Corollary~\ref{lagp}.
\hfill $\Box$
\vspace{.5cm}

Now let us relate this to $HPTS(d, A, n+1)_m$.

When~$G(d,A,n+1)_m < {{n+m} \choose m}$,
we have
\begin{equation}\label{Gfla}
G(d,A,n+1)_m=\sum_{t=1}^\infty (-1)^{t-1} \sum_{\ell(B)=t}  
{{n+|B|-t -(t-1)m}\choose n},
\end{equation}
where the inner sum is over subindices~$\overline{0} < B \subseteq A$.

We may  interpret each term as:

\medskip

\noindent {$t=1$: } $\displaystyle \sum_{i=1}^d {{n+k_i-1} \choose n}$, which is the
degree of the scheme;

\medskip

\noindent{$t=2$:} $\displaystyle \sum_{\{1 \leq i < j \leq d\}} {{n+k_i+k_j-2-m} \choose n}$,
counts (and subtracts) obstructions due to lines between pairs of points; 

\medskip

\noindent{$t=3$:}  $\displaystyle \sum_{\{1 \leq i_1  < i_2 < i_3 \leq d\}} 
{{n+k_{i_1}+k_{i_2}+k_{i_3} -3-2m} \choose n}$, counts obstructions
{\it to} lines between pairs of points given by planes between threesomes,
as a correction to the term $t=2$;

\noindent and so on.

That is, 
the $t$th term in~(\ref{Gfla}) accounts for each of the~$(t-1)$-planes 
occurring in the base locus of an~$A$-scheme in degree~$m$.

From the case of~$d \leq n+1$, where the SFI hypothesis does hold valid, e.g., 
we see that the function~$G$ exhibits the intersection numbers with regard to
the (multiple) planes determined by Corollary~\ref{lagp}.
Whence, we may view
the Strong Fr\"oberg-Iarrobino Conjecture as asserting that
a generic~$A$-scheme has  only linear
obstructions, and the Weak Conjecture as enumerating the linear obstructions occurring.
We shall exhibit this phenomenon in Section~6.

This suggests the following:

\begin{conj}\label{conj1} Let~$n, m, d \in \nats$, and~$A \in \nats^d$.

Suppose that~${\cal Z}$ is an~$A$-subscheme of~$\proj^n$.
Then
$$h_{\proj^n}({\cal Z},m) \leq G(d,A,n+1)_m.$$

Furthermore, assume that
$$G(d,A,n+1)_m < {{n+m} \choose m}.$$
We have equality of~$h_{\proj^n}({\cal Z},m)$ and $ G(d,A,n+1)_m$
precisely when~${\cal Z}$ exhibits only the expected linear obstructions
in degree~$m$.
\end{conj}

By definition,  if a multiple scheme has only the expected linear
obstructions in a given degree, none of its multiple subschemes can be nonlinearly
obstructed in this degree.   Under the viewpoint that $G(d,A,n+1)_m - HPTS(d,A,n+1)_m$
does count nonlinear obstructions to an $A$-scheme in degree~$m$ (as we shall justify in Section~6), let us observe:

\begin{cor}
Suppose that Conjecture~\ref{conj1} does hold valid.

Let~$n,m, d \in \nats$ and~$A \in \nats^d$.
Assume that~$G(d,A,n+1)_m < {{n+m} \choose m}$.

Suppose that ${\cal Z} \subset \proj^n$ is an~$A$-subscheme.
If 
$$h_{\proj^n}({\cal Z},m)= G(d, A,n+1)_m$$
then for every~$A' \leq A$ and each~$A'$-subscheme~${\cal Z}'$ of~${\cal Z}$ (respecting ordering)
we have
$$h_{\proj^n}({\cal Z}',m)=G(d, A',n+1)_m.$$
\end{cor}

This suggests, more generally, that the difference between the Hilbert function of a collection
of multiple points and the conjectural value keeps track of nonlinear obstructions, in the following sense:

\begin{conj}\label{conj2}  Let~$n,m, d \in \nats$ and~$A \in \nats^d$, for which
$G(d,A,n+1)_m < {{n+m} \choose m}$.  Take~$A' \in \nats^d$ for which~$A' \leq A$.

Suppose that~${\cal Z} \subset \proj^n$ is an~$A$-subscheme and that~${\cal Z}' \subset {\cal Z}$ 
is an~$A'$-subscheme.

Suppose that  for some~$\alpha \in \nats$ we have
$$h_{\proj^n}({\cal Z}',m) = G(d,A',n+1)_m - \alpha.$$  Then
$$h_{\proj^n}({\cal Z},m)  \leq  G(d,A,n+1)_m-\alpha.$$
\end{conj}

\medskip

Of course the conclusion of each conjecture is obvious in case the function~$G$ 
predicts maximal rank for a given scheme!
We ask the reader to check that this is not obvious in general!

Note also that Conjecture~\ref{conj2} is stronger than the 
Weak Fr\"oberg-Iarrobino Conjecture:  
given a generic $A$-scheme~${\cal Z}$
satisfying the hypothesis of the conjecture,  we do have 
$$h_{\proj^n}({\cal Z}_{red},m)=\deg {\cal Z},$$
so the conjecture predicts that $h_{\proj^n}( {\cal Z},m) \leq G(d,A,n+1)_m.$
However, experimental evidence on construction(!) points towards 
the statement's being simpler to verify inductively than the Weak Conjecture itself.

\section{Main Theorem}

The Main Theorem (Theorem~\ref{ubda}) gives an upper bound on the Hilbert function of any collection of
infinitesimal neighbourhoods of points in~$\proj^n$ based on Hilbert functions of certain such
subschemes of~$\proj^{n-1}$.
Particularly, the scheme of interest is shown in Theorem~\ref{siwa} to verify the bound given by
the Weak Fr\"oberg-Iarrobino Conjecture when each of the specified (``smaller'') ones  satisfy
the Strong Hypothesis.
Moreover, we obtain equality if and only if these schemes in lower dimension have only the expected
linear obstructions, as the conjectural function has been shown to compute.

The proof is attained, inductively, by the comparison in Lemma~\ref{key} of the Hilbert
function of a given~$A$-scheme with that of a $B$-scheme, where 
$B=A-(0, \ldots, 0, 1)$,  
To relate these we evaluate in Lemma~\ref{degree} the degree of the linear obstruction scheme
occurring between $A$ and~$B$, which may naturally be seen in terms
of a Hilbert function of fat points in codimension~one. 
Further, we see from Lemma~\ref{nlobs} that equality in the estimate of Lemma~\ref{key}
arises precisely when only the expected linear obstructions occur.

\begin{lemma}\label{degree}
Let $n, d, a \in \nats$, and $(j_1, \ldots, j_d) \in \nats^d$.
Let $p \in \aff^n$ and let  $L_1, \ldots, L_d$ be distinct lines
of~$\aff^n$ through~$p$.
Let~${\rho} \subset \aff^n$ be
the scheme $(p^a \cup L_1^{j_1} \cup \ldots \cup L_d^{j_d}) \cap p^{a+1}$.
In the projective space~$\proj^{n-1}$ of lines through~$p$ take the $(j_1, \ldots, j_d)$-scheme~${\cal W} \subset \proj^{n-1}$
given by~$L_1^{j_1}  \cup \ldots \cup L_d^{j_d}$.  Then $$\deg {\rho} = 
\deg p^a+h_{\proj^{n-1}}({\cal W},a).$$
\end{lemma}

{\it Proof:}  Call ${\frak m}$ the maximal ideal of $p$ in the
affine coordinate ring of $\aff^n$ and $I_1, \ldots, I_d$ the ideals
of the $d$ lines.  Identify ${\frak m}^a/ {\frak m}^{a+1}$ with the
vector space of forms of  degree $a$ in the prescribed projective space of lines
through $p$ and  and 
$({\frak m}^a \cap I_1^{j_1} \cap \ldots \cap  I_d^{j_d}+{\frak m}^{a+1})/{\frak m}^{a+1}$
 with forms of degree $a$ 
vanishing on the subscheme ${\cal W} \subset \proj^{n-1}$.
\hfill $\Box$

\medskip

\begin{lemma}\label{nlobs}
Let $n,m,k \in \nats$. 
Suppose that~${\cal Z} \subset \proj^n$ is any subscheme and~$p \in \proj^n$ so that~$p \notin {\cal Z}$.

Take~$\gamma$ as the intersection of~$p^k$ with the base locus of~$m$-ics through~${\cal Z} \cup \{p\}^{k-1}$
(so that $h_{\proj^n}({\cal Z} \cup \{p\}^{k-1},m)= h_{\proj^n}({\cal Z} \cup \gamma,m)$).
Then
$$ h_{\proj^n}( {\cal Z} \cup \{p\}^{k},m) 
=\min \left(h_{\proj^n}({\cal Z} \cup p^{k-1},m) + \deg p^k - \deg \gamma , {{n+m} \choose m}\right).$$
\end{lemma}

\medskip

\noindent{\bf Remark:}  
%
%
Quite generally, consider the base locus of (say) the linear system of~$m$-ics 
through a subscheme~${\cal Z} \subset \proj^n$.  
We may of course determine the base locus scheme~${\cal Y} \subset \proj^n$ from the~$m$th graded
piece~$I({\cal Z})_m$ of the ideal of~${\cal Z}$.
The subtlety is that there is no guarantee that~$I({\cal Y})_m = I({\cal Z})_m$.
Indeed, this issue may be viewed as the crux of the challenge in verifying the Segre conjecture
(along with analogues in higher dimension).

The point of Lemma~\ref{nlobs} is simply that we may locate such subschemes in a relative sense: 
comparing~$I({\cal Z} \cup \{p\}^{k-1})_m \supseteq I({\cal Z} \cup \{p\}^k)_m$ for~$p \in \proj^n$.
Whence we obtain a scheme~$\gamma$ with~$\{p\}^{k-1} \subset \gamma \subset \{p\}^k$
accounting for the base locus with respect to~$\{p\}^{k-1}$.  So, in the context of Lemma~\ref{key} we may
compare such a scheme~$\gamma$ with an expected linear obstruction scheme.

\medskip

{\it Proof:}  We may assume that 
$$h_{\proj^n} ({\cal Z} \cup \{p\}^k, m) < {{n+m} \choose n};$$
in particular, that $m \geq k$.

Choose co\"ordinates on~$\proj^n$ so that~$I:=I(p)=(X_1, \ldots, X_n)$.

Set $V=I({\cal Z} \cup p^{k-1})_m \cap I(p^k)_m$, so that
$ I^k_m \subseteq V \subseteq I^{k-1}_m. $

%
%
We may write
$$I^{k-1}_m = I^{k}_m+ X_0^{m-k+1} I^{k-1}_{k-1},$$
and hence we may find a linearly independent set~$F_1, \ldots, F_r \in (X_1, \ldots, X_n)^{k-1}$, each homogeneous of degree~$k-1$
so that
$$V=I^k_m+ X_0^{m-k+1} (F_1, \ldots, F_r).$$
So for the scheme $\gamma$ given by the ideal $(X_1, \ldots, X_n)^k+(F_1, \ldots, F_r)$, we have
$$I({\cal Z} \cup \gamma)_m = I({\cal Z} \cup \{p\}^{k-1})_m.$$ 
Note that~$\deg \gamma = \deg \{p\}^k - r$.

Corresponding to the forms~$F_1, \ldots, F_r$ we have forms~$G_1, \ldots, G_r \in I({\cal Z} \cup \{p\}^{k-1})_m$,
each distinguished by 
$G_i - X_{0}^{m+1-k} F_i \in I^k$.
In particular, no linear combination of the forms~$G_1, \ldots, G_r$ vanishes on~$\{p\}^k$, whence we do have
$$h_{\proj^n}({\cal Z} \cup \{p\}^k,m)= h_{\proj^n}({\cal Z} \cup \{p\}^{k-1},m) + r = h_{\proj^n}({\cal Z} \cup \{p\}^{k-1},m)+
\deg \{p\}^k - \deg \gamma.$$ 
\hfill $\Box$

We apply these results toward linear obstruction schemes:

\begin{lemma}\label{key}  Let~$n, m, d \in \nats$, and 
$A=(k_1, \ldots, k_{d}, k) \in \nats^{d+1}$.
Let~$\{p_1, \ldots, p_d, p\} \subset \proj^n$ and ${\cal Z} = \bigcup_{i=1}^d \{p_i\}^{k_i}$.
Choose a hyperplane $\proj^{n-1} \subset \proj^n$
for which~$\Gamma \cap \proj^{n-1} = \emptyset$. 
Take 
$$C=(c_1, \ldots, c_d) :=(k_1, \ldots, k_d)-\overline{m+1-k}.$$
and let $${\cal W}= \bigcup_{i=1}^d (\spax (p_i, p) \cap \proj^{n-1})^{c_i}.$$
Then
$$h_{\proj^n}( {\cal Z} \cup \{p\}^k, m) \leq h_{\proj^n}({\cal Z} \cup \{p\}^{k-1},m) + {{n+k-2} \choose {n-1}}
-h_{\proj^{n-1}}({\cal W} \cap \proj^{n-1}, k-1).$$   
Equality occurs precisely if ${\cal Z} \cup \{p\}^k$ is only linearly obstructed by~${\cal Z} \cup \{p\}^{k-1}$.

In particular,for a generic scheme,  take~$B=(k_1, \ldots, k_d, k-1)$. Then:
\begin{eqnarray}
HPTS(d+1,A,n+1)_{m} &\leq& 
HPTS(d+1,B,n+1)_{m}+{{n+k-2} \choose {n-1}}-
\nonumber \\
&&HPTS(d,C,n )_{k-1}.
\nonumber
\end{eqnarray}
Equality occurs exactly when a generic~$A$-subscheme of~$\proj^n$ 
is only linearly obstructed by a~$B$-subscheme.

\end{lemma}

{\it Proof:}
Let us put together the relevant information from our previous observations.

Take $\gamma$ as the intersection of~$\{p\}^k$ with the
base locus of the linear system of $m$-ics through~${\cal Z} \cup \{p\}^{k-1}$.
So~$\gamma$ contains the linear obstruction scheme~$\rho$ which is
(by Lemma~\ref{lobs}) given by
$$\rho= (p^{k-1} \cup {\cal W}) \cap p^k.$$

By Lemma~\ref{nlobs} we have

$$h_{\proj^n}({\cal Z} \cup \{p\}^k,m ) 
= h_{\proj^n}( {\cal Z} \cup \{p\}^{k-1},m) + \deg \{p\}^k - \deg \gamma,$$
so that 
$$h_{\proj^n}({\cal Z} \cup \{p\}^k,m ) 
\geq h_{\proj^n}( {\cal Z} \cup \{p\}^{k-1},m) + \deg \{p\}^k - \deg \rho,$$
with equality occurring exactly when~${\cal Z} \cup \{p\}^k$ is only linearly obstructed by~${\cal Z} \cup \{p\}^{k-1}$.

From Lemma~\ref{degree} we may now plug in~$\deg \rho = \deg \{p\}^{k-1} + h_{\proj^{n-1}}({\cal W} \cap \proj^{n-1},k-1)$
to obtain the desired conclusion.
\hfill $\Box$   

\medskip

Altogether we have:

{\it Proof of Theorem~\ref{ubda}:}

Let $n, m, d \in \nats$ and~$A=(k_1, \ldots, k_d) \in \nats^d$.  
Let us take~$C_{ji}$ as stated; that is, $C_{ji} = (k_1, \ldots, k_j)+\overline{i-m}$, for~$j=1, \ldots, d-1$,
and~$i=0, \ldots, k_j-1$.

Given an~$A$-subscheme~${\cal Z} \subset \proj^n$ we wish to produce~$C_{ji}$-subschemes~${\cal W}_{ji} \subset \proj^{n-1}$
from~${\cal Z}$ for which
\begin{equation} \label{leq}
h_{\proj^n}({\cal Z}, m) \leq \deg {\cal Z}- \sum_{j=1}^{d-1} \sum_{i=1}^{k_j -1} h_{\proj^{n-1}} ({\cal W}_{ji}, i);
\end{equation}
equality holding just when~${\cal Z}$ displays only the expected linear obstructions in degree~$m$.
Namely, taking~$\proj^{n-1} \subset \proj^n$ as a  hyperplane that does not meet the support of~${\cal Z}$
we have
$${\cal W}_{ji} = (\bigcup_{r=1}^j \spax (p_{j+1},p_r)^{k_r+i-m} ) \cap \proj^{n-1}.$$
From Lemma~\ref{key} we obtain this by double induction on~$d$ (starting with~$d=0$) and then on~$k_d$
(from the initial value~$k_d=0$).

We find that equality holds in~(\ref{leq}) exactly when~${\cal Z}$ exhibits only the expected linear obstructions,
again from Lemma~\ref{key}.
\hfill $\Box$

\medskip

Let us compare the description in Theorem~\ref{ubda} of the Hilbert function of an~$A$-subscheme in~$\proj^n$
with the behaviour of the conjectural function of Fr\"oberg:

\begin{lemma}\label{Gcomp}
 Let $n, m, d \in \nats$, $A=(k_1, \ldots, k_d)$, and~$C_{ji}=(k_1, \ldots, k_j)+\overline{i-m}$,
for~$j=1, \ldots, d-1$, $i=0,\ldots, k_j-1$.
If~$G(d, A, n+1)-m < {{n+m} \choose m}$ then
$$G(d, A, n+1)_m = \sum_{j=1}^{d-1} \sum_{i=1}^{k_j-1} G(j,C_{ji},n)_i.$$
(More technically, it is enough to assume that~$G'(\ell(A'), A', n+1) \leq {{n+m} \choose m}$ for
each~$A' \leq A$.)
\end{lemma}

{\it Proof:} Compute directly from equation~(\ref{Gfla}). 
\hfill $\Box$

\medskip

Hence, according to Theorem~\ref{ubda} the function~$G$ may be viewed directly  as keeping track
of the expected linear obstruction schemes identified in Definition~\ref{exp}.  
We find, then:

{\it Proof of Theorem~\ref{siwa}:}
Let us take~$n,m, d \in \nats$, $A=(k_1, \ldots, k_d) \in \nats^d$, and~$C_{ji}=(k_1, \ldots, k_j)+\overline{i-m}$
for each~$j=1, \ldots, d-1$ and~$i=0, \ldots, k_j-1$.

Suppose, as indicated, that
$$HPTS(j,C_{ji},n)_i = G(j, C_{ji},n)_i,$$
for each~$j=1, \ldots, d-1$ and~$i=0,, \ldots, k_j-1$.

If~$G(d, A, n+1)_m = {{n+m} \choose n}$ then the Weak Conjecture holds trivally for an~$A$-scheme.
Otherwise, take~${\cal Z} \subset \proj^n$ as a generic~$A$-subscheme.  By Theorem~\ref{ubda}
 we have
\begin{eqnarray}
h_{\proj^n}({\cal Z},m) &\leq&\sum_{j=1}^d {{n+k_j-1} \choose n} - \sum \sum G(j,C_{ji},n)_i, 
\nonumber
\\
&=& G(d,A,n+1)_m,
\end{eqnarray}
and equality holds exactly when~${\cal Z}$ has only the expected linear obstructions in degree~$m$.
  Hence, by upper-semicontinuity, the inequality
applies to every~$A$-subscheme of~$\proj^n$ and the Weak Conjecture is satisfied by each~$A$-scheme in degree~$m$.

Let us now take~${\cal Z}_0$ as an arbitrary~$A$-subscheme.
 of~$\proj^n$, and take~${\cal W}_{ji}$ as
the~${\cal C}_{ji}$-subschemes of~$\proj^{n-1}$ identified in Theorem~\ref{ubda}.
Note that~$h_{\proj^{n-1}}({\cal W}_{ji},i) \leq G(j,C_{ji},n)_i$ for each pair~$j,i$.

 Suppose that~${\cal Z}_0$ has only the expected linear obstructions in degree~$d$.
Then $$h_{\proj^n}({\cal Z},m) = \deg Z - \sum \sum h_{\proj^{n-1}} ({\cal W}_{ji},i) 
\geq G(d,A,n+1)_m.$$
According to the Weak Conjecture, as we've just verified in this case, we must have equality here.
\hfill$\Box$


\medskip

\vspace{.2cm}

Now let us examine the special case of quasihomogeneous schemes (including homogeneous schemes, 
as in the Strong Conjecture):

\medskip
\begin{prop}\label{siww}
Let $n \in \nats$.
Suppose that  the Strong Fr\"oberg-Iarrobino Conjecture holds in $\proj^{n-1}$.
Assume further that the Weak Fr\"oberg-Iarrobino Conjecture applies
to each  quasihomogeneous
scheme of fat points in $\proj^n$ with support on $n+4$ points.
Then \\
a) the  Weak Fr\"oberg-Iarrobino Conjecture holds valid in $\proj^{n}$ for every
 quasihomogeneous fat point scheme in~$\proj^n$, and \\
b) a generic quasihomogeneous
collection of infinitesimal neighbourhoods of points in~$\proj^n$ satisfies the
Strong Hypothesis if and only if it exhibits only the expected linear obstructions in
each degree~(Definition~\ref{exp}).\\
c) for any quasihomogeneous collection of infinitesimal neighbourhoods in~$\proj^n$, its
Hilbert function agrees with the function given by the
Strong Hypothesis provided that it has only the expected linear obstructions.
\end{prop}

Hence to verify the Weak Conjecture in this setting, along with examining the Strong Conjecture itself,
 it suffices to examine schemes
supported on~$n+4$ points of~$\proj^n$. 

\medskip 

{\it Proof:}
We verify in Theorem~\ref{p3} (Section~7) that the Weak Conjecture holds valid (indeed, for schemes of mixed multiplicities)
in~$\proj^n$ for each~$n \leq 3$.  Further, according to results from~\cite{mernc} along with the given hypothesis,
we have that each fat point subscheme of~$\proj^n$ supported on at most~$n+3$ points satisfies
the Weak Conjectures.

Note, next,  that in order to apply Theorem~\ref{ubda} to a quasihomogenous subscheme of~$\proj^n$, 
it suffices to see that the Strong Conjecture (on homogeneous schemes) applies to  
corresponding subschemes of~$\proj^{n-1}$.

According to the given assumptions  along with the conditions given by the Strong Conjecture,
it remains to deal with~$\proj^n$ for~$n =4$ and $n= 5$.

Let us make the following observations:
\begin{itemize}
\item Fix~$n,m, \ell \in \nats$.  The function~$G(d, (\ell, \ell, \ldots, \ell, k),n+1)_m$ is strictly
increasing in~$k$ until it reaches its maximum value~$\displaystyle{{n+m} \choose m}$.
\item Consider the ``extra'' exceptions predicted by the Strong Conjecture occurring for~$n \leq 4$
(given by homogeneous schemes supported on at least~$n+5$ points). 
For~$n=3$, the additional cases (i.e., with~$d \geq 8$) have~$d=9$ and  $m=2k$,
and the extra cases for~$n=4$ ($d \geq 9$)
are: $d=14$, $m=2k$, $k=2$ or~$3$.  In each of these, we have
$$G(d,k,n+1)_{m} = {{n+m} \choose m}.$$ 
\end{itemize}

Hence, each time such a case arises in the application of Theorem~\ref{siwa} toward a 
quasihomogeneous $A$-scheme, in a given degree~$m$, we have already 
that~$$G(d, A, n+1)_m = {{n+m} \choose m}.$$

More precisely, take~$n=4$.  According to Theorem~\ref{siwa} we must show that
$$HPTS(10, (\ell, \ldots, \ell, k),5)_m \leq G(10, (\ell, \ldots, \ell, k),5)_m,$$
whenever~$2(\ell+k-1-m)=k-1.$  
Note that for each~$r \leq k-2$ we have~$2(\ell+r -m) < r$, so we do have
$$HPTS(10, (\ell, \ldots, \ell, k-1),5)_m \leq G(10, (\ell, \ldots, \ell, k-1)_m.$$
Further, by the previous observations,
$$G(10, (\ell, \ldots, \ell, k-1)_m = {{m+4} \choose 4},$$
and hence the same holds for~$G(10,(\ell, \ldots, \ell, k)_m$, and we are done.

Similarly for~$n=5$. 
\hfill$\Box$.

\medskip


\medskip
Let us remark on the comparison between determining when the Strong Hypothesis applies to
an arbitrary scheme and finding obstructions to that scheme:

\begin{prop}\label{new}
Suppose  that the Weak Fr\"oberg-Iarrobino Conjecture holds valid; that is, the Hilbert function of 
any collection of multiple points is bounded above by the conjectural value.

Let $n,m, d \in \nats$ and~$A \in \nats^d$.
Suppose that~${\cal Z} \subset \proj^n$ is any~$A$-subscheme.
Then~${\cal Z}$ satisfies the Strong Fr\"oberg-Iarrobino Hypothesis (resp., in a given degree)
provided that~${\cal Z}$ presents only the expected (linear) obstructions (resp., in that degree).
\end{prop}

{\it Proof:} 
It suffices to prove the result in a given degree~$m$.

We may assume without loss of generality 
that there is an~$m$-ic 
vanishing on~${\cal Z}$ and that~$n \geq 2$.

Take~${\cal Z} \subset \proj^n$ as prescribed and~${\cal W}_{j,i} \subset \proj^{n-1}$
as in Theorem~\ref{ubda}.  
So~${\cal Z}$ has only the expected linear obstructions in degree~$m$ if and only if
$$h_{\proj^n}({\cal Z},m)= \deg {\cal Z} - \sum \sum h_{\proj^{n-1}}({\cal W}_{ji},i).$$
Now according to the Weak Conjecture, applied to~$\proj^{n-1}$, we have that
\begin{equation}\label{impt}
\sum \sum h_{\proj^{n-1}} ({\cal W}_{ji},i)  \leq \sum \sum G(j, C_{ji}, n)_i
\end{equation}
Hence (from the Weak Conjecture applied to~$\proj^n$ along with Lemma~\ref{Gcomp})
if~${\cal Z}$ has only the expected linear obstructions in degree~$m$ we have equality
in~(\ref{impt}), so that
$$h_{\proj^n}({\cal Z},m)= G(d,A,n+1)_m.$$
Notice, it follows that if~${\cal Z}$ has the only expected linear obstructions, 
so must each of the schemes~${\cal W}_{ji}$.
\hfill $\Box$

\medskip

Let us compare our results here with the conjectures overall.
Take an $A$-subscheme of~$\proj^n$.
Then Theorem~\ref{ubda} provides the conclusion of the Weak Conjecture when each
 of the derived~$C_{ji}$-schemes in~$(n-1)$-space do satisfy
the Strong Fr\"oberg-Iarrobino Hypothesis.
Otherwise,  each~$C_{ji}$-scheme
suffering from lack of $SFI$ must have nonlinear obstructions.
Conjecture~\ref{conj1} demands that each such obstruction is then carried over to the~$A$-scheme.


%
\section {The Weak Conjecture holds valid in~$\proj^3$}

We illustrate here the use of Theorem~\ref{ubda} in verifying the Weak Fr\"oberg-Iarrobino
Conjecture in a given dimension without the full requirement of the Strong Hypothesis in
each lower dimensional case.

We prove Theorem~\ref{p3}, that
the Weak-Fr\"oberg conjecture does hold in~$\proj^n$ for~$n \leq 3$.
Of course, since the Strong Conjecture of Fr\"oberg and Iarrobino
holds in~$\proj^1$, regardless of multiplicities, then 
Weak Conjecture does in~$\proj^2$.
However, as described in Section~3, the Strong Conjecture remains open in~$\proj^2$ and 
presents many exceptions to the Strong Hypothesis.  
We observe, though,  that
in the homogeneous situation
the required results for the application
of Theorem~\ref{ubda} follow easily from general results of
Nagata on multiple points in~$\proj^2$.  So, to deal with mixed multiplicities
we make use of the following:

\begin{lemma}(Numerical Observation)
 Let $A=(k_1, \ldots, k_d) \in \nats^d$, with $k_1 \geq k_2+2$.
Take $B=(k_1-1, k_2+1, k_3, \ldots, k_d)$.
Then $G(d, A, n+1)_m \geq G(d,B,n+1)_m$.

Indeed, equality  occurs if and only if $G(d,B,n+1)_m = {{n+m} \choose m}$.
\end{lemma}

(For the main idea, note that for a fixed pair of integers $n, c$ the
maximal value of the quantity $\displaystyle {{n+a} \choose n}+{{n+c-a} \choose n}$  
is obtained from $\displaystyle a=\left\lceil \frac c 2 \right\rceil$, so $\displaystyle c-a=\left\lfloor \frac c 2 \right\rfloor$.) 

\medskip

\medskip

{\it Proof of Theorem~\ref{p3}:} 
Let~$A=(k_1, \ldots, k_d, \ell) \in \nats^{d+1}$, with~$k_1 \geq k_2 \geq \ldots \geq k_d \geq \ell$.
Fix~$m \in \nats$.  We wish to show that 
$$HPTS(d+1, A, 4)_m \leq G(d+1, A, 4)_m.$$
Let us assume inductively that such an inequality holds for each~$B < A$; 
without loss of generality we may also assume that~$G(d, A, 4)_m < {{m+3} \choose 3}$.

First, if~$k_1+\ell \leq m+1$ we may take~$B=(k_1, \ldots, k_d)$ and then 
$$G(d+1, A, 4)_m  = G(d,B,4)_m+ {{\ell+2} \choose 3}, $$
and we are done by the induction hypothesis with respect to~$B$.

So let us assume that~$k_1 + \ell \geq m+2$.
We now take~$B=(k_1, \ldots, k_d, \ell-1)$ and then~$C=(c_1, \ldots, c_s) \in \nats^s$ by rewriting 
the~$d$-uple~$(k_1, \ldots, k_d)-\overline{m+1-\ell}$ so that only positive terms occur.
(Namely, we 
take~$s \leq d$ maximal with respect to the property that~$k_s+k \geq m+2$.
and then~$C=(k_1, \ldots, k_s)-\overline{m+1-k}$. ) 
By Lemma~\ref{key} we have
$$HPTS(d+1,A,4)_m \leq HPTS(d+1,B,4)_m + {{\ell+1} \choose 2} -HPTS(s, C, 3)_{\ell-1}.$$
By the induction hypothesis, we are done once we see that~$HPTS(s, C, 3)_{\ell-1}=G(s, C, 3)_{\ell-1}$.

If~$s \leq 4$ we are done. (See Section~4.)

Next, let us observe that~$s \leq 6$.    Note, first, that  for each~$a \in \nats$ and each~$m \leq 2a-1$, 
\begin{equation}\label{e3}
G(8,\overline{a}, 4)_m = {{m+3} \choose 3}. 
\end{equation}
(To see this it suffices to evaluate~$G(8,\overline{a},4)_m$ 
for~$m=2a-1$,
where~$G$ computes the degree of the scheme.)
%
Hence if~$ s \geq 7$, we would have
$$G(d, A, 4)_m \geq G(8,(k, \ldots, k, \ell),4)_m,$$ where~$k=k_s$.
It is easy to see, by the numerical observation lemma, that we may find an integer~$a$
for which the latter item is at least~$G(8,\overline a,4)_m$ and~$2a \geq m+1$.
So  equation~(\ref{e3}) contradicts our hypothesis on~$G(d,A,4)_m$.

Thus we are left with the cases~$5 \leq s \leq 6$.  
From Nagata's results~\cite{n} (see Section~3) we have 
$$HPTS(s, C, 3)_{\ell-1} =G(s,C,3)_{\ell-1}$$ provided that:
$$\sum_{i=1}^5 (k_i+\ell-1-m) \leq 2(\ell-1)+1.$$

Let us claim, then, that we do have this inequality due to the hypotheses on~$A$.
That is, imagine that $\sum_{i=1}^5 (k_i+\ell-1-m) \geq 2(\ell-1)+2.$
Let us check that $$G(d+1, A, 4)_m = {{m+3} \choose 3}.$$
Since~$A \geq (k_1, \ldots, k_5, \ell)$ it is enough to replace~$A$ by the latter.

As before we may find an integer~$a$ for which $G(d,A,4)_m \geq G(6,\overline a, 4)_m$
and $5 (2a-1-m) \geq 2(a-1)-3$, so $8a \geq 5m$.  In particular $2a \geq m+2$ 
but $3a < 2m$ and we have
$$G(6,\overline a, 4)_m=6 {{a+2} \choose 3} - { 6 \choose 2} {{2a+1-m } \choose 3}.$$
It is easy to compute (substitute  $\displaystyle m=\lfloor{\frac {8a} 5 }\rfloor $, say)
 that the latter quantity is at least $\displaystyle {{m+3} \choose 3}$, as claimed.
\hfill $\Box$

\section{General consequences}

Here we shall derive some immediate consequences of the results obtained in Section~6 toward~$\proj^n$.
We start with a focus on the  ``first order'' cases of the Weak Fr\"oberg-Iarrobino conjecture;
namely, where expected
base loci do not include planes of dimension two.
By means of Macaulay duality these  give a direct analogue to the result of Iarrobino 
(Theorem~\ref{ith})
on generic forms.  
For this, we compare (via Lemma~\ref{lconj2}) the behaviour
of the Hilbert function in a given degree~$m$ with that of each~$(A-\overline i)$-scheme
in degree~$(m-i)$, respectively.
This provides a verification of Conjecture~\ref{conj1} in the first-order situation.

Let us start in the following four observations by applying Theorem~\ref{siwa} directly. 

\begin{cor} Let~$A=(k_1, \ldots, k_d) \in \nats^d$.
Assume that $k_1 \geq k_2 \geq \ldots \geq k_d$.  
Then
$$HPTS(d,A,n+1)_m \leq
G(d,A,n+1)_m $$
for $m \geq k_1+k_2-3$.
\end{cor}

{\it Proof:} If $G(d, A,n+1)_m = {{n+m} \choose m}$
(such as in the case~$m \leq k_1-1$) we are done.

Assume, by induction on the number of points and the orders,
 that the result holds for all $B < A$ (such as $|B|=1$).
 
Call~$B=(k_1-1, k_2, \ldots, k_d)$,
and~$C=(c_2, \ldots, c_d)$ where~$c_i=k_1+k_i-1 -m$.
So~$c_i \leq 2$ for~$i=2, \ldots, d$.
 
We have  
$$G(d,A,n+1)_m = G(d,B,n+1)_m+{{n+k-2} \choose {n-1}}-G(d-1,C,n)_{k-1},$$
and $G(d,A,n+1)_m > G(d,B,n+1)_m$,
so that $G(d-1,C,n)_{k-1} < {{n+k-2} \choose {n-1}}$.
By the Alexander-Hirschowitz theorem~(\cite{h},\cite{a}, \cite{ah1}, \cite{ah2}, \cite{ah3})
we have $$HPTS(d-1,C,n)_{k-1}=G(d-1,C,n)_{k-1}.$$
 
Whence by Lemma~\ref{key} we find
$$HPTS(d,A,n+1)_m \leq HPTS(d,B,n+1)_m + {{n+k-2} \choose {n-1}}-
G(d-1,C,n)_{k-1},$$
which by the induction hypothesis is at most:
$$G(d,B,n+1)_m +{{n+k-2} \choose {n-1}}-G(d-1,C,n)_{k-1}=G(d,A,n+1)_m.$$
\hfill $\Box$


\begin{cor} Let~$n, d, k \in \nats$.
Then
$$HPTS(d,\overline{k},n+1)_{2k-2} \leq {{n+k-1}\choose n}d-{d \choose 2}.$$
\end{cor}

\begin{cor}  Let~$n, d, k \in \nats$. Assume that $k \geq 4$.  Then
$$HPTS(d,\overline{k},n+1)_{2k-3} \leq {{n+k-1}\choose n}d- (n+1){d \choose 2}.$$
\end{cor}

\begin{cor}  Let~$n, d, \in \nats$ and $A=(k_1, \ldots, k_d) \in \nats^d$.
 If $\max\{k_i\} \leq 4$ then 
for each degree~$m$ we have 
$$HPTS(d,A,n+1)_m \leq G(d,A,n+1)_m.$$
\end{cor}

\begin{cor} Let~$n, d, k \in \nats$.

Assume that for each~$d_1 \leq d$, $k_1 \leq k-1$, and~$m_1 \geq 2 k_1 -1$ 
 a generic union of~$d_1$
$k_1$-uple points of~$\proj^{n-1}$ has maximal rank with respect to
$|{\cal O}_{\proj^{n-1}}(m_1)|$ for every degree $m_1 \geq 2k_1-1$.

Then in each degree $m$ with $ 2m \geq 3k-2 $ we have
$$HPTS(d,k,n+1)_{m} \leq G(d,k,n+1)_{m}.$$
\end{cor}

Likewise one obtains an analogous conclusion on an $A$-scheme of mixed 
multiplicities, $A=(k_1, \ldots, k_d)$ where $k_1 \geq \ldots \geq k_d$
and the degree~$m$ under consideration satisfies~$2m \geq k_1+2k_2 -2$.
We shall see in  Corollary~\ref{nohyp} how to simplify the hypotheses of the
above result and apply toward $A$-schemes, homogeneous or otherwise.

Let us continue the examination of situations in which the base locus of a 
system of $m$-ics through a general collection of fat points is expected to
contain lines but not planes.  We aim toward simplifying the use of Theorem~\ref{ubda}
under such a circumstance.
In Corollary~\ref{pith} this gives a result directly comparable to that of Iarrobino
in the setting of the 
Fr\"oberg conjectures.

The main instrument is the following:

\begin{lemma}\label{lconj2}
Let $n, m, d \in \nats$ and~$A \in \ints^d$.
Given a generic~$A$-subscheme~${\cal Z} \subset \proj^n$, let~${\cal Z}^{-r}$
denote the corresponding~$(A-\bar r)$-subscheme of~${\cal Z}$ for
each~$r \in \nats$.

If 
$$h_{\proj^n}({\cal Z}^{-1}, m-1) = \deg {\cal Z}^{-1} - \alpha$$
then
$$h_{\proj^n}({\cal Z},m) \leq \deg {\cal Z}-\alpha.$$

Therefore, if
$$h_{\proj^n}({\cal Z},m) = \deg {\cal Z}$$
we have
$$h_{\proj^n}({\cal Z}^{-r},m-r) = \deg {\cal Z}^{-r}$$
for each~$r=0, \ldots, m$.

%
%
\end{lemma}

{\it Proof:} 
Let $\Gamma = {\cal Z}_{red}$, and write $A=(k_1, \ldots, k_d)$.  
We may assume without loss of generality that~$k_i \geq 2$ for each~$i$ from~1 to~$d$.

Assume that~$h_{\proj^n}({\cal Z}^{-1},m-1) = \deg {\cal Z}^{-1} - \alpha$.

Take the homogeneous co\"ordinate ring~$S={\cal K}[X_0, \ldots, X_n]$ of~$\proj^n$ with 
co\"ordinates chosen so that none of the points~$[1:0:\ldots:0]$, up to~$[0:\ldots :0:1]$
lie on~${\cal Z}$.
View~$R={\cal K}[X_1, \ldots, X_n]$ as the co\"ordinate ring of the hyperplane~$\proj^{n-1}$ described
by the form~$X_0$, and take~$\pi_p$ as the projection from the points~$p=[1:0:\ldots:0]$ onto~$\proj^{n-1}$.

We obtain then that the ideal~$I(\pi_p(\Gamma)) = I(\Gamma) \cap R$ has $I(\Gamma) \cap R_m = d$.
Consider the exact sequences in the commutative diagram:

$$
\begin{CD}
\\ 0 @>>>
R_m  @>>>
S_m
@>{\displaystyle {{\frac {\partial} {\partial X_0} }}}>>
S_{m-1} @>>> 0  \\
& &\bigcup
 & &\bigcup & & ||  & \\
0 @>>>  I(\Gamma)_m \cap R_m @>>>   I(\Gamma)_m @>{\displaystyle {\frac {\partial} {\partial X_0} }}>>  
S_{m-1} @>>> 0 \\
\end{CD}
$$

\medskip

Let us filter:
$$S_m \supseteq  I(\Gamma)_m= V_0 \supset V_1 \supset \ldots \supset V_n = I(\z)_m,$$

where
$$\displaystyle V_j= \left\{ F \in V_0 \; : \; {\frac {\partial F} {\partial X_i} } \in  I(\z^{-1}),
 i=1, \ldots, j \right\}.$$

\medskip

From the diagram:

$$
\begin{CD}
\\ 0 @>>>
V_0 \cap R_m  @>>>
V_0
@>{\displaystyle {{\frac {\partial} {\partial X_0} }}}>>
S_{m-1} @>>> 0  \\
& &\bigcup
 & &\bigcup & & \bigcup  & \\
0 @>>>  V_1 \cap R_m @>>>   V_1 @>{\displaystyle {\frac {\partial} {\partial X_0} }}>>  
I({\cal Z}^{-1})_m  @>>> 0 \\
\end{CD}
$$

\medskip

\noindent we see that~$\dim V_0/V_1 \leq \deg \z^{-1} - \alpha$.
Routinely we obtain
$$\dim V_j /V_{j+1} \leq  \sum_{i=1}^d {{n+k_i-1-j} \choose {n-j}} $$ 
for~$j=1, \ldots, n-1$ (namely, the degree of an~$(A-\overline{1})$-subscheme of~$\proj^{n-j}$.
 In sum we then obtain the desired inequality.
\hfill $\Box$

\medskip

\noindent{\bf Remark:} The argument applies equally to an~$A$-scheme with arbitrary support~$\Gamma \subset \proj^n$
provided that the projection~$\pi_p(\Gamma) \subset \proj^{n-1}$ does have~$h_{\proj^{n-1}} (\pi_p(\Gamma),m)=d$.

\begin{cor}\label{nohyp}
Let~$n, m, d \in \nats$, and~$A = (k_1, \ldots, k_d) \in \nats^{d}$.
Suppose that
$$k_1 \leq \ldots \leq k_d \mbox{~and~} 2 k_{d-1}+k_d \leq 2m+2.$$
Let
$$C=(c_1, \ldots, c_{d-1}):=(k_1, \ldots, k_{d-1})-\overline{m+1-k_d}.$$
Suppose that
$$HPTS(d, C,n)_{k_d-1}=G(d,C, n)_{k_d-1}.$$
Then
$$HPTS(d+1,A,n+1)_m \leq G(d+1,A,n+1)_m.$$
Again, equality applies to a generic~$A$-scheme 
if and only if no obstructions occur other than the expected linear ones.
\end{cor}

\noindent{\bf Remark.}  We may likewise extend the corollary to an arbitrary~$A$-subscheme of~$\proj^n$
under the hypothesis that the scheme~${\cal W} \subset \proj^n$
identified in Lemma~\ref{key} achieves the value of~$G$ in degree~$k_{d}-1$, and the projection of~${\cal Z}_{red}$
to~$\proj^{n-1}$ attains maximal rank in degree~$m$.

\medskip

{\it Proof:}  We may assume that~$G(d, A, n+1)_m < {{n+m} \choose m}$.

Take
$C_{j,i}=(k_1, \ldots, k_{j-1})-\overline{m-i}$, for each~$j=2, \ldots, d$
and~$i=0, \ldots, k_{j}-1$.

By Theorem~\ref{ubda} we are done once we see that
$$HPTS(j, C_{j,i}, n)_{i} = G( j, C_{j,i}, n)_{i}$$ for all~$j,i$.

By hypothesis, equality holds for $(j,i)=(d,k_d-1)$.
So by Lemma~\ref{lconj2} we obtain equality as well
in the case of $(d-1,i)$ for all~$i \leq k_d-1$.
 
According to our numerical hypothesis
this says that for each~$i \leq k_d-1$,
a~$C_{d,i}$-scheme imposes independent conditions
on~$i$-ics.
Hence its subschemes, notably, 
the $C_{j,i}$-schemes for~$j=1, \ldots, k_d-1$ do as well.

\hfill $\Box$

\noindent{\bf Remark.}  Now let us compare the result of Corollary~\ref{nohyp}
with Theorem~\ref{ith} of Iarrobino (see Section~4).

Take $$A=(k_1, \ldots, k_d, k), k \geq \max{k_i},$$ and
$$C=(k_1, \ldots, k_d)-\overline{m+1-k},$$ 
as in the statement of Corollary~\ref{nohyp}.

Then the Macaulay dual of~$A$ in degree~$m$ is described by
the~$(d+1)$-uple
$$A^{\perp} = (j_1, \ldots, j_d, j)$$
where  $j_i=m-k_i+1$
and~$j=m+1-k$, so~$j \leq \min\{j_i\}$. 
The dual of~$C$ in degree~$k-1$ ($=m-j$) is given by
$$C^{\perp}=(j_1, \ldots, j_d).$$

The hypothesis that~$2m \geq k_{d-1}+k_d+k-2$ 
(along with~$m \leq k_d+k-1$, to make things interesting, say)
gives that~$j_d+j \leq m \leq j_{d-1}+j_{d}+j-1$, i.e., we are
in the range described by Iarrobino.

From his statement 
that
$$HGEN(d+1, A^{\perp},n+1)_m \geq F(d+1,A^{\perp},n+1)_m$$
if 
$$HGEN(d,C^{\perp},n+1)_{m-j} = F(d,C^{\perp},n+1)_{m-j}.$$
one may expect to obtain 
the upper bound~$HPTS(d+1,A,n+1)_m \leq G(d+1,A,n+1)_m$ 
from information on a~$C$-scheme living in~$\proj^n$.
By Corollary~\ref{nohyp} we obtain:

\begin{cor} \label{pith} Let $n, m, d, k \in \nats$.  
Suppose that~$A \in \nats^{d+1}$ and $C \in \ints^d$ 
satisfy the numerical hypotheses of Corollary~\ref{nohyp}.

Assume that~${\cal Z} \subset \proj^n$ is a~$C$-subscheme 
supported on a generic subset of a hyperplane~$\proj^{n-1}$.
If
$$h_{\proj^n}({\cal Z},k-1)=G(d,C,n+1)_{k-1}$$ then
$$HPTS(d+1,A,n+1)_m \leq G(d+1,A,n+1)_m.$$
\end{cor}

{\it Proof:} One should only  notice from  the sequence
$$ 0 \rightarrow {\cal I}_{\tilde {\cal Z}}(k-2) \rightarrow
{\cal I}_{\cal Z}(k-1) \rightarrow 
{\cal I}_{{\cal Z} \cap \proj^{n-1},\proj^{n-1}} (k-1)
\rightarrow 0$$
we have~$H^1(\proj^{n-1},{\cal I}_{{\cal Z} \cap \proj^{n-1},\proj^{n-1}} (k-1)) = 0$,
so that~$HPTS(d,C,n)_{k-1} = G(d,C,n)_{k-1}$.

By Lemma~\ref{lconj2} we obtain that
$$HPTS(d,C-\bar i,n)_{k-1-i}=G(d,C-\bar i,n)_{k-1-i}$$ for each~$i$
from~$0$ to~$k-1$,
so 
$$ \sum_{i=0}^{k-1} H(d, C-\bar i,n)_{k-1-i} = G(d,C,n+1)_{k-1}.$$

Taking $B=(k_1, \ldots, k_d)$ we find
that 
$$HPTS(d+1,A,n+1)_m \leq HPTS(d,B,n+1)_m+{{n+k-1} \choose n} -G(d,C,n+1)_{k-1}.$$

Assuming inductively that
$HPTS(d,B,n+1)_m \leq G(d,B,n+1)_m$ we see that, as advertised,
$$HPTS(d+1,A,n+1)_m \leq G(d+1,A,n+1)_m.$$
\hfill$\Box$



\section{Verifying the Strong Fr\"oberg-Iarrobino Hypothesis}

We illustrate here examples of the agreement of the Hilbert function of a generic collection of
fat points  with the value predicted by the Strong Fr\"oberg-Iarrobino
Hypothesis, using Castelnuovo lower bounds.

In Proposition~\ref{rnc} we see that the Hilbert function of a collection
of fat points has the predicted value in a given degree when the sum of
its multiplicities is ``not too large'' compared with the degree and the dimension
of the projective space..
(Indeed, the numerical hypotheses are geared so that,
 according to conjectures of \cite{cg}, every $A$-scheme
with support in linearly general position should have this property.)

Next, in Proposition~\ref{plus1} we reexamine cases, studied in~\cite{me}, of a homogeneous union of fat points of
multiplity~$k$ in degree~$k+1$, the least degree in which the Hilbert function is 
not obviously attainable.
From lower bounds in~\cite{me} together with upper bounds given by Theorem~\ref{ubda}
we find equality in the cases addressed.

To start, let us recall that 
if~${\cal Z} = {\cal Y} \cup {\cal X}$, 
in which~${\cal X}$ is an~$A$-subscheme of~$\proj^n$
 supported on~$\proj^{n-1}$
and no component of~$\cal Y$ has support on~$\proj^{n-1}$ then 
\begin{equation}\label{cast}
h_{\proj^n}({\cal Y} \cup {\cal X}, m) 
\geq h_{\proj^n}({\cal Y} \cup \tilde{\cal X}, m-1)
+h_{\proj^{n-1}}({ \cal X} \cap \proj^{n-1}, m),
\end{equation}

where~$\tilde{\cal X}$ is the~$(A-\bar 1)$-subscheme of~$\proj^n$
given by~${\cal I}_{\tilde{\cal X}} = {\cal I}_X: {\cal I}_{\proj^{n-1}}$.

We shall refer to this as a {\bf Castelnuovo lower bound}.

\begin{prop}\label{rnc}
Let $n, d, m \in \nats$.
Let $k_1, \ldots, k_d$ be nonnegative integers so that
$k_1+\ldots+k_d \leq mn+1$ (or $d \leq n+1$).  
Then 
$$HPTS(d, (k_1, \ldots, k_d), n+1)_m =G(d,(k_1, \ldots, k_d),n+1)_m.$$
\end{prop}

{\it Proof:} Let us order the orders so that~$k_1 \geq k_2 \geq \ldots \geq k_d$.
Consider an~$A$-subscheme of~$\proj^n$ with~$A=(k_1, \ldots, k_d)$.
Let us assume by induction that the conclusion of the theorem holds:
\begin{itemize}
\item in $\proj^{n-1}$, inducting on~$n$ (as in~$n=1$), and
\item in degree~$m-1$, inducting on degree (such as $m \leq k_1-1$).
\item in each case of a~$B$-subscheme of $\proj^n$ having~$B < A$, 
by induction on $\max \{k_j\}$ and then by induction on~$d$
(such as  $d =1$).
\end{itemize}

We  assume, then, that~$d \geq 2$,
and~$m \geq k_1$.  
Further, according to Corollary~\ref{lagp} 
we may assume that $\sum_{i=1}^d k_i \leq mn$ in case~$d \leq n+1$. 

Let us start with the case that~$m \leq k_1+k_2-1$.  
(The case~$m \geq k_1+k_2-1$ has already been proved by Catalisano,
Trung, and Valla in \cite{ctv}; but we shall include a proof for
completeness.)

Assume $\sum_{i=1}^d k_i \leq mn+1$.

We start by showing that $HPTS(d, A, n+1)_m$ is bounded below by~$G(d,A,n+1)_m$,
by applying the Castelnuovo inequality~(\ref{cast}) 
toward the inductive hypotheses.
To do this it is necessary to distinguish between the cases $k_1=m$ and~$k_1 < m$, 
in order to ensure that the values of $G$  obtained in the lower dimension and degree
do add up to $G(d,A,n+1)_m$.

\noindent{\bf Case:} $k_1=m$.

Take a general~$A$-scheme~${\cal Z}=\{p\}^m \cup {\cal Z}_1$.
By Lemma~\ref{lag} the (scheme-theoretic) cone ${\cal C}$ over~${\cal Z}$ through~$\{p\}^m$
is in the base locus of the system of~$m$-ics through~${\cal Z}$, and (so)
$$h_{\proj^n}({\cal Z},m)=h_{\proj^n}({\cal C},m).$$
Choose a general hyperplane~$\proj^{n-1} \subset \proj^n$ (namely, not through~$p$!).
Since the restriction map
$${\cal I}_{\cal C}(m) \rightarrow {\cal I}_{{\cal C} \cap \proj^{n-1}}(m) \rightarrow 0$$
is surjective on global sections we have
\begin{eqnarray}
h_{\proj^n}({\cal Z},m)&=& h_{\proj^n}( {\cal C},m)
\nonumber 
\\
&=&
 h_{\proj^n}({\cal C},m-1) + h_{\proj^{n-1}}({\cal C} \cap \proj^{n-1},m).
\nonumber 
\end{eqnarray}

Now ${\cal C} \cap \proj^{n-1}$ is a $(k_2, \ldots, k_d)$-subscheme of~$\proj^{n-1}$
for which~$k_2+\ldots+k_d \leq m(n-1)+1$, while
$$h_{\proj^n}({\cal C},m-1)={{n+m-1} \choose m}=G(d, A, n+1)_{m+1}$$
so that by induction on~$n$ (along with inspection on the behaviour of~$G$) we obtain
$$HPTS(d,A,n+1)_m = G(d,A,n+1)_m.$$

\noindent{\bf Case:} $k_1 \leq m-1$:

\noindent{\bf Claim.}  Under the given hypotheses,
suppose that $\z=\z_0 \cup \z_1$, where $\z_1$ is a generic
$(k_3, \ldots, k_d)$-subscheme of $\proj^n$ supported on~$\proj^{n-1}$
and $\z_0$ is a general
$(k_1, k_2)$-subscheme of  $\proj^n$.
Then
$$h_{\proj^n}(\z,m) = G(d,A,n+1)_m.$$
 
To see this, fix $k_1 \geq k_2$.
Assume by induction that the claim holds valid in degree~$m-1$.

Let $L=\spax\{p_1,p_2\}$ and $\{q\} = L \cap \proj^{n-1}$, so
by Lemma~\ref{lag} we have
$$h(\z, m)= h(\z \cup L^{k_1+k_2-m},m).$$
According to Castelnuovo, we have
$$h(\z, m) \geq  h(\z_0 \cup \tilde \z_1,m-1) +
h((q^{k_1+k_2-m} \cup \z_1) \cap \proj^{n-1},m).$$
By our hypothesis on the claim, we have
$$h(\z_0 \cup \tilde \z_1,m-1) \geq
G(d, (k_1, k_2, k_3-1, \ldots, k_d-1),n+1)_{m-1}$$
since $\sum_{i=1}^d k_i -(d-2) \leq (m-1)n+1$. 
 
Further, by induction on $n$ we have
$$h((\{q\}^{k_1+k_2-m} \cup \z_1) \cap \proj^{n-1},m)=
G(d-1, (k_1+k_2-m, k_3, \ldots, k_d), n)_m,$$
since
$(k_1+k_2-m)+k_3+\ldots+k_d \leq m(n-1)+1$.
 
Hence
\begin{eqnarray}
h(\z,m) &\geq&
G(d, (k_1, k_2, k_3-1, \ldots, k_d-1),n+1)_{m-1}+
\nonumber \\
&&~
G(d-1, (k_1+k_2-m, k_3, \ldots, k_d), n)_m
\nonumber \\
&=& G(d, A, n+1)_m.
\nonumber
\end{eqnarray}

Whence, by upper semicontinuity we have 
$HPTS(d,A,n+1)_m \geq G(d,A,n+1)_m$.

Now to see that equality holds, take~$B=A-(1, 0, \ldots, 0)$
and~$C=(k_2, \ldots, k_d)-\overline{m+1-k_1}$.
Rewrite~$C$ as~$C=(c_2, \ldots, c_r)$, where~$r $ is 
maximal for which~$k_1+k_r > m+1$.
One may easily compute that $|C| \leq (k_1-1)(n-1)$ so that
by induction on~$n$ we have $$HPTS(r-1,C,n)_{k_1-1}=G(r-1,C,n)_{k_1-1},$$
and hence we may apply Lemma~\ref{key} to obtain the inequality on $HPTS(d,A,n+1)_m$.

%
%

Let us consider the case $m \geq k_1+k_2$, as in \cite{ctv}.
Here~$G(d,A,n+1)_m$ is the degree of an~$A$-scheme, hence
$HPTS(d,A,n+1) \leq G(d,A,n+1)_m$.
Given a generic $A$-scheme~$\z$ take a hyperplane~$H$ containing
exactly~$n$ of the reduced points of~$\z$.
Since~$\sum_{i=1}^d k_i \leq mn+1$ we have $\sum_{i=1}^{d-n} (k_i-1)+\sum_{i=d-n+1}^d k_i \leq (m-1)n+1$ 
so 
$$h(\z_0 \cup \tilde \z_1, m-1) =
G(d, (k_1, \ldots, k_{d-n}, k_{d-n+1} -1, \ldots, k_d-1),n+1)_{m-1}.$$
We certainly have
$$h(\z_1 \cap H,m)=G(n,(k_{d-n+1}, \ldots, k_d),n)_m,$$
so that~$h({\cal Z}_0 \cup {\cal Z}_1,m)=G(d,A,n+1)= \deg \z$.

 \hfill $\Box$

\noindent{\bf Example.} Consider the exceptional cases to the Strong
Fr\"oberg-Iarrobino Conjecture, given by a~$\overline{k}$-scheme supported on~$d=n+3$ points
of~$\proj^n$, in a given degree~$m$.
According to Proposition~\ref{rnc}, we do have equality between
$HPTS(n+3,A,n+1)_m$ and~$G(n+3,A,n+1)_m$ unless~$|A| \geq mn+2$.
But in that case the scheme does have a subscheme, given by intersection
with the rational normal curve passing through, which does not
present maximal rank in degree~$m$, and yet meets linear obstruction
subschemes transversely.
One should expect then that this does give an exception to the
strong conjecture.

We observe the following (which subsumes Corollary~\ref{nplus3}):

\begin{cor}\label{n++3} 
Let~$n, m \in \nats$, and~$A \in \nats^{n+3}$.
Suppose~$G(n+3,A,n+1)_m < {{n+m} \choose m}.$
 
We have~$HPTS(n+3, A, n+1)_m = G(n+3,A,n+1)_m$ 
provided that~$|A| \leq mn+1$.

Indeed, suppose that that Conjecture~\ref{conj1} holds. 
Then~$HPTS(n+3, A, n+1)_m = G(n+3,A,n+1)_m$ if and only
if~$|A| \leq mn+1$.

\end{cor}

We shall see in~\cite{mernc} that the conclusion of the corollary does hold valid (regardless of assumptions
on Conjecture~\ref{conj1}) by closely examining fat points on a rational normal curve.

\medskip

Now let us move on to  verifying the conjecture in the setting
of points of equal multiplicity in the lowest nontrivial degree.

\begin{cor}\label{plus1} Let~$n, k \in \nats$.
Suppose $\displaystyle{ d \leq \max \left(n+1, {\frac {(n+3)(n+2)} { 2(k^2-1)} }  \right).}$
Then $$HPTS(d,\overline{k},n+1)_{k+1}= G(d,\overline{k},n+1)_{k+1}.$$
\end{cor}

{\it Proof:} We have seen in \cite{me} that 
\begin{equation}\label{up}
HPTS(d,\overline{k},n+1)_{k+1} \geq G(d,\overline{k},n+1)_{k+1}.
\end{equation}

If~$d \leq n+1$ we are done.

Otherwise, for each~$j \leq k-1$ we have
$${\frac { (n+3)(n+2) } { 2(k^2-1)} }  < 
{\frac { (n+2)(n+1)}  { 2((k-j)^2-1)}} ,$$
so that by~(\ref{up}) we have that
$$HPTS(d-1,\overline{j-1},n)_j \geq G(d-1,\overline{j-1},n)_j.$$
Therefore by Theorem~\ref{siwa}
we obtain the inequality
$$HPTS(d, \overline{k},n+1)_{k+1} \leq G(d, \overline{k},n+1)_{k+1},$$
and hence equality.
\hfill $\Box$

\medskip

\noindent{\bf Remark.}  Notice that in degree~$k+1$ it should suffice to deal with~$\overline{k}$-schemes supported
on~$d$ points, with
$${{n+k-1} \choose n} d \leq {{n+k+1} \choose n};$$
that is,
$$d \leq {\frac {(n+k+1)(n+k)} {(k+1)k } }$$
so that Corollary~\ref{plus1} covers ``about half the ground''.

\section{Toward Refinement of the Strong Fr\"oberg-Iarrobino Conjecture}
We consider here the issue of when the Strong Hypothesis should hold valid.
We construct counterexamples to the Strong Conjecture in~$\proj^n$, for~$n= 4, 5$, and~$6$.
Namely, for each~$n$ we exhibit~$k(n)$ so that for each~$k \geq k(n)$ a generic collection of~$n+5$ $k$-uple points does not satisfy
the conjecture.

In the homogeneous situation, the main cases neglected by the SFI conjecture in~$\proj^n$
are given on~$n+3$ or~$n+4$ points.
For mixed multiplicities, Corollary~\ref{nplus3} says that the hypothesis does apply to~$n+3$ fat points
in degree~$m$ provided that the scheme meets the rational normal curve through these
points to degree at most the Hilbert function in that degree.
Otherwise, we do expect that the curve presents nonlinear obstructions causing the
failure of the SFI hypothesis.
Likewise for~$n+4$ fat points we should find that ``most'' exceptions to the SFI hypothesis
are due to excess intersection with a rational normal curve or an elliptic normal curve.
But notice, numerically, that the latter already implies the former.
(So that the next main suspect would perhaps be a singular curve.)

Let us determine conditions under which 
a generic homogeneous scheme of~$d$ fat points in~$\proj^n$ should not satisfy the Strong Fr\"oberg-Iarrobino Conjecture.
Consider the following inequalities:

\begin{equation} \label{max}
G(d,\overline{k}, n+1)_m < {{n+m} \choose n}
\end{equation}
\begin{equation} \label{rn1}
(n+3)k \geq mn+2,
\end{equation}
and
\begin{equation}\label{rn2}
(d-1)(2k-1-m) \leq (k-1)(n-1)+1.
\end{equation}

Suppose that a triple~$(n,k,m)$ satisfies each of (\ref{max}), (\ref{rn1}), and (\ref{rn2}) for a given number of points, $d$.
Let us observe that a generic~$\overline{k}$-subscheme on~$d$ points of~$\proj^n$ disobeys the 
Strong Fr\"oberg-Iarrobino Hypothesis under such conditions.

From~(\ref{rn2}), along with Proposition~\ref{rnc} we have that
$$HPTS(j,\overline{k+i-m},n)_i = G(j, \overline{k+i-m},n)_i$$
for each~$j=1, \ldots, d-1$ and~$i=0,\ldots,k-1$. 
By~(\ref{rn2}), a generic~$\overline{k}$-subscheme on~$d$ points of~$\proj^n$ is nonlinearly obstructed,
so from Theorem~\ref{siwa} along with~(\ref{max}) we have that
$$HPTS(d, \overline{k}, n+1)_m < G(d, \overline{k}, n+1)_m,$$
as asserted.

Let us review the Strong Fr\"oberg Conjecture (Conjecture~\ref{spt}). 
The assertion is that each homogeneous collection of~$d$ multiple points in~$\proj^n$ has Hilbert
function that agrees with the corresponding function~$G$, except when: $d=n+3$ or~$n+5$; or else~$n \leq 4$ when
there are futher exceptions.  For~$n=4$ the additional exceptions are given when~$d=9$ (multiplicity either~2 or~3). 
                                                                 
We now exhibit counterexamples to the conjecture, in each dimension~${n=4, 5, 6}$ and for~$d=n+5$.
In each of these cases, the Strong Conjecture predicts that a generic homogeneous~$A$-scheme has Hilbert function
equal to the function~$G$.
For each~$n, k$ take~$m(n,k)$ as the greatest  integer~$m$ for which: $mn \leq (n+3)k-2$, 
so that the triple~$(n,k,m(n,k))$ satisfies inequality~(\ref{rn1}).

Notice that for~$n=4$ or~$5$, each triple~$(n,k,m(n,k))$ also satisfies~(\ref{rn2}), 
while~$(6,k,m(6,k))$ satisfies~(\ref{rn2}) provided that~$k$ is odd.

Next, define~$k(n)$ as the minimal integer so that: for each~$k \geq k(n)$ we have:
\begin{equation}
(n+5) {{n+k-1} \choose n} \leq {{n+m(k,n)} \choose n},
\end{equation}
(provided that such an integer exists).
If so, inequality~(\ref{max}) applies 
to each triple~$(n,k,m(n,k))$ when~$k \geq k(n)$. 

We may compute: $k(4)=88, k(5)=88$, $k(6)= 141$;  $k(7)$ is between 231 and 648.
(One may check, as well, that many counterexamples occur for~$k < k(n)$.)

To summarise:

\begin{prop} \label{ctr}
If~$n= 4, 5, $ or~$6$ and~$k \geq 88,\; 88$,  or $ 141$, respectively, a generic collection of~$n+5$ $k$-uple points in~$\proj^n$
violates the Strong Fr\"oberg Conjecture.
\end{prop}  

{\it Proof:} We have seen the cases of $n=4$, $n=5$, (and also $n=6$ in the case of~$k$ odd).
However, Proposition~\ref{rnc} does not guarantee the remaining cases.  
We shall nonetheless deduce the required conclusion using the proposition. 

To complete the case of~$n=6$, it remains to obtain the following:

\medskip

\noindent{\bf Claim:}  A generic union of~10 $r$-uple points in~$\proj^5$ imposes independent conditions in 
each degree~$ \geq 2r-1$ for~$r \geq 3$,
that is,
$$HPTS(10, \overline{r},6)_{2r-1} = 10 {{r+4} \choose 5}.$$

So after verifying the claim, it follows from Theorem~\ref{siwa} that, for each~$k \geq 141$ we have:
$$HPTS(10, \overline{k}, 7)_{m(k,n)} < G(10,\overline{k},7)_{m(k,6)}.$$

To verify the claim, let us choose a flag~$\proj^3 \subset \proj^4 \subset \proj^5$, along with a
general hyperplane~${H \subset \proj^5}$, and take~$\Gamma = \Gamma_1 \cup \Gamma_2$
 so that~$\Gamma_2$ consists of~$5$ points
on~$\proj^3$ and~$\Gamma_1$ contains~5 points on~$H$, chosen sufficiently generally.
So we make the further:

\medskip

\noindent{\bf Claim:}  For each $r \geq 3$ we have:
 $$h_{\proj^5}(\Gamma_1^r \cup \Gamma_2^r,2r-1)= 10  {{r+4} \choose 5}.$$

To see this, we start by taking~${\cal S}$ as the union of lines between points of~$\Gamma_1$.
By~(\ref{cast}), we have
$$h_{\proj^5} (\Gamma_1^r \cup \Gamma_2^r, 2r-1) \geq 
h_{\proj^5}(\Gamma_1^{r-1} \cup \Gamma_2^{r-1}, 2(r-1)-1)+
h_H(\Gamma_1^r, 2r-2)+ h_{\proj^4}( {\cal S} \cap \proj^4 \cup \Gamma_2^r, 2r-1).$$

As we have seen (e.g, Section~5) we obtain that 
$$h_H(\Gamma_1^r, 2r-2)= 5 {{3+r} \choose 4} - {{5} \choose 2},$$
for each~$r \geq 2$.

Next, by~(\ref{cast}), 
$$h_{\proj^4}( \proj^4 \cap {\cal S} \cup \Gamma_2^r, 2r-1) \geq 
\sum_{j=0}^{r-1} h_{\proj^3}(\Gamma_2^{r-j}, 2r-j-1)+
h_{\proj^4}(\proj^4 \cap {\cal S} ,r-1).
$$

To each item in the sum, Proposition~\ref{rnc} applies ($r \geq 2$), and we have
$$h_{\proj^4}({\cal S} \cap \proj^4, r-1) ={5 \choose 2},$$
when~$r \geq 3$.

It remains then to verify the initial case~$r=3$ of the claim.
For this, let us specialise $\Gamma_1= \Sigma \cup \{p\}$ so that~$ p \in H \cap \proj^4$
and take~${\cal S}_1$ as the union of lines through points of~$\Sigma$.
Now:
\begin{equation}
 h_{\proj^5}(\Gamma_1^3 \cup \Gamma_2^3, 5) \geq  
h_{\proj^5}( \Sigma^3 \cup \{p\}^2 \cup \Gamma_2^2, 4)+
h_{\proj^4}(\proj^4 \cap {\cal S}_1 \cup \{p\}^3 \cup \Gamma_2^2,4)+
h_{\proj^3}(\Gamma_2^3,5).
\nonumber
\end{equation}
The last term in the sum displays (as we have observed) independent conditions, equal to~$50$.

For the middle term, consider a further (but temporary) specialisation of~$p$ to a generic point~$q$
of~$H \cap \proj^3$ so that the term is bounded below by the sum:
$$h_{\proj^4}(\proj^4 \cap {\cal S}_1 \cup \{q\}, 2)+ h_{\proj^3}(\{q\}^2 \cup \Gamma_2, 3)+
h_{\proj^3}(\{q\}^3 \cup \Gamma_2^2,4),$$
each term of which contributes the expected number of conditions, so that we do have
$$h_{\proj^4}(\proj^4 \cap {\cal S}_1 \cup \{p\}^3 \cup \Gamma_2^2, 4)= 46.$$

We are left with
$$h_{\proj^5}(\Sigma^3 \cup\{p\}^2 \cup \Gamma_2^2, 4) \geq 
h_{\proj^5}(\Sigma^2 \cup \{p\} \cup \Gamma_2^2, 3) + h_{H}(\Sigma^3 \cup \{p\}^2, 4),$$
where~Proposition~\ref{rnc}  applies to the latter term.  For the former, we observe that
the technique of \cite{me3} applies directly.  We obtain, then,
$$h_{\proj^5}(\Sigma^3 \cup \{p\}^2 \cup \Gamma_2^2, 4) = 55+59=114, $$
as desired.

The grand total is then:
$$h_{\proj^5}(\Sigma^3 \cup \{p\}^3 \cup \Gamma_2^3, 5) \geq 114+46+50=210=10 {7 \choose 2},$$
so that equality holds and the claim is verified.
\hfill$\Box$  

\medskip

\noindent{\bf Remark:}  Let us now look at the case~$n=7$, where~$k(7) \leq  648$, and for each~$k \geq k(n)$
we have that~$(n,k,m(n,k))$ satisfies both~(\ref{max}) and~(\ref{rn1});  however~(\ref{rn2}) shall not apply.
So we do not immediately obtain counterexamples from Proposition~\ref{rnc}, and carrying out further
investigation seems a little less easy than for~$\proj^6$.

But note: in each case a generic~$\overline{k}$-scheme on~$12$ points of~$\proj^7$ exhibits nonlinear
obstructions in degree~$m(n,k)$.  Hence, each case must either contradict the Strong Fr\"oberg Conjecture
or Conjecture~\ref{conj1}!

(We suggest to the curious reader to extend the argument of Proposition~\ref{ctr} to showing that~11 
generic~$\overline{k}$-uple points in~$\proj^6$ do entertain the SFI Hypothesis in each degree~$m$ 
for which~$4m \geq 7k$ as follows: Choose a flag~$\proj^3 \subset \proj^4 \subset \proj^5 \subset \proj^6$.
Then apply the technique (repeatedly) to a 
collection~$\Sigma=\Gamma_6 \cup \Gamma_5 \cup \Gamma_4 \cup \Gamma_3$ of eleven points, 
where~$\Gamma_j \subset \proj^j$, for each~$j$, $|\Gamma_6| = |\Gamma_5|=|\Gamma_4|=2$,
and~$\Gamma_3|=5$.  
Notice that Proposition~\ref{rnc} does apply to~$\proj^3 \cap \Gamma_3^k$ in degree~$m$.
Check then whether (or when) the method of Proposition~\ref{rnc} does apply by splitting each
restriction into sums.) 

\medskip

In general we seek to identify conditions under which the Strong Hypothesis should apply to collections of points 
mixed multiplicities, as is done in~\cite{cm1} in the case of~$\proj^2$.
Geometrically, then, this requires the study of the base locus of a linear system with multiple base points;
specifically its excess intersection with a (nonlinear) positive dimensional variety.
To start, results of~\cite{mernc} yield the Hilbert functions of collections (of arbitrary multiplicities) lying
on a rational normal curve.  
One obtains predictions of when a  a collection of multiple points of~$\proj^n$ disobeys the Strong Hypothesis
on account of a rational normal curve through~$n+3$ of the reduced points.
The main principle is then to identify subvarieties of~$\proj^n$ for which excess intersection with a general 
collection of infinitesimal neighbourhoods impedes the Strong Hypothesis on its Hilbert function, as in Conjecture~\ref{hh}.

\section{The algebraic reinterpretation}
Let us look back to the algebraic conjectures of Fr\"oberg and
Iarrobino.  Particularly, we shall stick to the cases of appropriate
characteristic and number of points.

As remarked earlier, the Strong Fr\"oberg Conjecture derives
from the expectation that the minimal free resolution of an
ideal~$I$ generated by general forms should exhibit only Koszul
relations  ``as much as possible'' with respect to degrees.
But, of course, $\dim I_m$ needn't itself predict the resolution
in degree~$m$.

However, in the case of the Strong Algebraic Fr\"oberg-Iarrobino conjecture
we obtain the extra information from Theorem~\ref{ubda}.
Namely, when we dualise to the study of multiple points, the
equality of~$HPTS$ with the function~$G$ should predict
(according to Conjecture~\ref{conj1}) only linear
obstructions.  Redualising and interpreting linear obtstructions
by Macaulay methods, such an equality yields the information
that the corresponding syzygies are exactly as predicted.
Whence the Strong Conjecture of Fr\"oberg-Iarrobino not only implies that
of Fr\"oberg but gives desired conclusions on the resolution.

\begin{prop}  Assume that Conjecture~\ref{conj1} holds.

Suppose that~$I \subset S={\cal K}[X_0,\ldots, X_n]$
is an ideal generated by powers of linear forms, and  that~$I$ satisfies
the Strong Algebraic Fr\"oberg-Iarrobino Conjecture.

Take $M$ maximal for which~$I_M \ne S_M$ (as determined by the conjecture).
Then for each~$m \leq M$ the $m$th graded piece of the mnimal free
resolution is Koszul.
\end{prop}

Notice how Proposition~10.1 compares with a recent result of~\cite{mmr}.
There an ideal $(F_1, \ldots, F_d)$ given by a general collection of
forms is considered under the hypothesis that each of the ideals~$(F_1, \ldots, F_j)$
satisfies the Strong Fr\"oberg Conjecture for~$j=1, \ldots, d$.
The conclusion on resolution is then just as in the proposition above.

So our assumptions in the proposition are partly stronger, partly weaker 
than those in~\cite{mmr}, but with the bonus of geometric insight.

\medskip
\medskip
\small{
\address{
DEPARTMENT OF MATHEMATICS, UNIVERSITY OF NOTRE DAME, NOTRE DAME, IN 46556\\
{\it Email address:} kchandle@@noether.math.nd.edu}}


\begin{thebibliography}{9999999}

\bibitem[A] {a} J. Alexander,
{\em Singularit\'es imposables en position g\'en\'erale aux
hypersurfaces de~$\proj^n$}, Compositio Math. {\bf 68} (1988), no. 3, 305-354.

\bibitem[AH1]{ah1} J. Alexander, A. Hirschowitz,
{\em Un lemme d'Horace diff\'erentiel: application aux
singularit\'es hyperquartiques de~$\proj^5$},
J. Alg. Geom.~{\bf 1} (1992), no.~3, 411-426.
\bibitem[AH2]{ah2} J. Alexander, A. Hirschowitz,
{\em La m\'ethode d'Horace \'eclat\'ee: application \`a
l'interpolation en degr\'e quatre},
Invent. Math. {\bf 107}, no, 3, 585-602 (1992).
\bibitem[AH3]{ah3} J. Alexander, A. Hirschowitz,
{\em Polynomial interpolation in several variables},
J. Alg. Geom. {\bf 4} (1995), no. 2, 201-222.
\bibitem[AH4]{ah5} J. Alexander, A. Hirschowitz,
{\em An asymptotic vanishing theorem for generic unions of multiple points},
Invent. Math {\bf 140} (2000), no. 2, 303-325.
\bibitem[An]{an} D. Anick, 
{\em Thin algebras of embedding dimension three}, 
J. Algebra {\bf 100} (1986), no. 1, 235-259.
\bibitem[Au]{au} M. Aubry, 
{\em S\'erie de Hilbert d'une alg\`ebre de polyn\^omes quotient},  J. Algebra {\bf 176} (1995), no. 2, 392-416. 
\bibitem[CEG]{cg} M.V. Catalisano, Ph. Ellia, A. Gimigliano, {\em Fat points on rational normal curves},
J. Algebra {\bf 216} (1999), no. 2, 600-619. 
\bibitem[CTV] {ctv} M.V. Catalisano, N.V. Trung, G. Valla, 
{\it A sharp bound for the regularity index of fat points in general position}, 
Proc. Amer. Math. Soc. {\bf 118} (1993), no. 3, 717-724. 
\bibitem[C1]{me0} K. Chandler, {\it Geometry of Dots and Ropes},
Trans. Amer. Math. Soc., {\bf 347 } (1995), no. 3, 767-784.
\bibitem[C2]{me} K. Chandler, {\it Higher infinitesimal neighbourhoods}, 
J. Algebra {\bf 205} (1998), no. 2, 460-479. 
\bibitem[C3]{me2} K. Chandler, 
{\it A brief proof of a maximal rank theorem for generic double points in projective space}, 
Trans. Amer. Math. Soc, {\bf 353} (2000), no. 5, 1907-1920. 
\bibitem[C4]{me3} K. Chandler, {\it Linear systems of cubics singular at general points of projective space}, 
Compositio Math., {\bf 134} (2002),  no. 3, 269-282.
\bibitem[C5]{me4} K. Chandler, {\it Uniqueness of exceptional singular quartic}, Proc. Amer. Math. Soc., to appear.
\bibitem[C6]{mernc} K. Chandler, {\it Multiple points on a rational
normal curve}, in preparation.
\bibitem[C7]{methree} K. Chandler, {\it Triple points in projective space}, in preparation.
\bibitem[Ci] {cil} C. Ciliberto, {\it Geometric aspects of polynomial interpolation in more variables and of Waring's problem},
 European Congress of
 Mathematics, Vol. I (Barcelona, 2000), 289-316, Progr. Math., 201, Birkhäuser, Basel, 2001.
\bibitem[CM1]{cm1} C. Ciliberto, R. Miranda,
{\it Degenerations of planar linear systems}, J. Reine Angew. Math 501
(1998), 191-220.
\bibitem[CM2]{cm2} C. Ciliberto, R. Miranda, 
{\it Linear systems of plane curves with base points of equal multiplicity},
 Trans. Amer. Math. Soc. 352 (2000), no. 9, 4037-4050.
\bibitem[CM3]{cms} C. Ciliberto, R. Miranda, {\it The Segre and Harbourne-Hirschowitz Conjectures}, 
Applications of algebraic geometry to coding theory, physics and computation (Eilat, 2001), 37-51, 
NATO Sci. Ser. II Math. Phys. Chem., 36, 
Kluwer Acad. Publ., Dordrecht, 2001. 
\bibitem[CCMO]{cmo} C. Ciliberto, F. Cioffi, R. Miranda, F. Orecchia,
{\it Bivariate Hermite interpolation via computer algebra and algebraic geometry techniques},
World Scientific Publ., Lecture Notes on Computing.
\bibitem[EI]{ei} J. Emsalem, A. Iarrobino, {\em Inverse system of a symbolic power. I},
J. Algebra {\bf 174} (1995), no. 3, 1080-1090. 
\bibitem[F]{f} R. Fr\"oberg,
{\it An inequality for Hilbert series of graded algebras},              
Math. Scand. 56 (1985), no. 2, 117-144. 
\bibitem[FH]{fh} R. Fr\"oberg, J. Hollmann,
{\em Hilbert series for ideals generated by generic forms}, 
J. Symbolic Comput. {\bf 17} (1994), no. 2, 149-157. 
\bibitem[GS]{gs} M. Gasca, T. Sauer, {\em Polynomial interpolation in several variables},
Advances in Symbolic Computation, {\bf 12} (2000) 377-410.
\bibitem[Ha]{h1} B. Harbourne, {\it Points in good position in $\proj^2$},
Zero-dimensional schemes (Ravello, 1992) 213-229, de Gruyter, Berlin, 1994.
\bibitem[Ha2]{h3} B. Harbourne, 
{\it An algorithm for fat points on $\bold P\sp 2$},
Canad. J. Math. 52 (2000), no. 1, 123-140.
\bibitem[Ha3]{h2} B. Harbourne, {\it Problems and Progress: A survey on fat points in $\proj^2$} 
Zero-dimensional schemes and applications (Naples, 2000), 85-132, 
Queen's Papers in Pure and Appl. Math., 123, 
Queen's Univ., Kingston, ON, 2002.
\bibitem[H]{h} A. Hirschowitz,
{\em La m\'ethode d'Horace pour l'interpolation \`a plusieurs variables},
Manuscripta Math. {\bf 50} (1985), 337-388.
\bibitem[HL]{hl} M. Hochster, D. Laksov,
{\it The linear syzygies of generic forms}, 
Comm. Algebra {\bf 15} (1987), no. 1-2, 227--239. 
\bibitem[I]{iar} A. Iarrobino, {\it Inverse system of a symbolic power III:
thin algebras and fat points}, Compositio Math. {\bf 108} (1997), no. 3, 319-356.
\bibitem[IK]{ik} A. Iarrobino, V. Kanev, {\it Power Sums, Gorenstein algebras, and determinantal loci},
Lecture Notes in Mathematics, 1721, Springer-Verlag, Berlin, 1999.
\bibitem[MM-R]{mmr} J. Migliore, R. Miro-R\'oig, {\it Ideals of general forms and the ubiquity of the
Weak Lefschetz Property}, to appear in Journal of Pure and Applied Algebra.
\bibitem[N]{n} M. Nagata, {\it On rational surfaces, II}, 
Mem. Coll. Sci. Univ. Kyoto. Ser. A. Math. {\bf 33} (1960), 271-293
\bibitem[S]{se} B. Segre, {\it Alcune questioni su insiemi finiti di punti in geometria algebrica }, Atti Convegno
intern. di Geom. Alg. di Torino, 15-33.
\bibitem[St]{s} R. Stanley, {\it Hilbert functions of graded algebras}. Advances in Math. {\bf 28} (1978), no. 1, 57-83.
\end{thebibliography}
\end{document}